\documentclass [12pt]{article}
\usepackage{amsmath,amssymb}

\begin{document}

   \title {Cramer's rules for   Hermitian systems  of coquaternionic 
 equations.}

\author {Ivan Kyrchei \footnote{Pidstrygach Institute
for Applied Problems of Mechanics and Mathematics, str.Naukova 3b,
Lviv, Ukraine, 79060, kyrchei@lms.lviv.ua}}
\date{}
\maketitle
\begin{abstract}
In this paper  properties of the determinant of a Hermitian matrix are investigated, and   determinantal representations of the inverse of a Hermitian coquaternionic matrix are  given. By their using, Cramer's rules for left and right systems of linear equations with Hermitian coquaternionic matrices of
coefficients are obtained. Cramer's rule for a two-sided coquaternionic matrix equation ${\bf AXB}={\bf D}$ (with Hermitian ${\bf A}$, ${\bf B}$) is given as well.
\end{abstract}

\section{Introduction}

 \newtheorem{thm}{Theorem}[section]
 \newtheorem{cor}[thm]{Corollary}
 \newtheorem{lem}[thm]{Lemma}
 \newtheorem{prop}[thm]{Proposition}

 \newtheorem{defn}[thm]{Definition}

 \newtheorem{rem}[thm]{Remark}
 \newtheorem{ex}{Example}
 \numberwithin{equation}{section}
A quaternion algebra ${\bf{H}}( a,b)$ over a field ${\bf{F}}$ (denoted by $(\frac{a,b}{{\bf{F}}})$) are a central simple algebra over ${\bf{F}}$, and a four-dimensional vector space over
${\bf{F}}$  with basis $\{1,i,j,k\}$
and the following multiplication rules:
\begin{gather*}
  i^{2}=a, \,\,\,
  j^{2}=b, \,\,\,
  ij=k, \,\,\,
  ji=-k,
\end{gather*}
where $\{a,b\}\subset {\bf{F}}$.
 To every quaternion algebra ${\bf{H}}( a,b)$, one can associate a quadratic form ${\bf n}$
(called the norm form) on ${\bf{H}}$ such that ${\bf
n}(xy)={ \bf{n}}(x){\bf n}(y)$, for all $x$ and $y$ in
${\bf{H}}$. A linear mapping $x \rightarrow \overline{x}={\bf
t}(x)-x$ is also defined on ${\bf{H}}$. It is an involution with properties
\begin{gather*}\overline{\overline{x}}=x,\, \overline {x + y} = \overline
{x} + \overline {y},\, \overline {x \cdot y} = \overline {y}
\cdot \overline {x}.\end{gather*}
 An element
$\overline{x}$ is called the conjugate of $x \in{\bf{H}}$. ${\bf t}(x)$ and ${\bf n}(x)$ are called the trace
and the norm of $x$ respectively. $\{{\bf n}(x), {\bf
t}(x)\} \subset {\bf{F}}$ for all $x$  in ${\bf{H}}$ and possess  the following conditions,
\begin{gather*}{\bf n}\left(
\overline{x} \right) = {\bf n}(x),\,{\bf t}\left( \overline{x}
\right) = {\bf t}(x),\, {\bf t}\left( {q \cdot p} \right) =
{\bf t}\left( {p \cdot q} \right).\end{gather*}
Depending on the choice of ${\bf{F}}$, $ a$ and $b$ we
have only two possibilities (\cite{le}):
\begin{enumerate}
  \item $(\frac{a,b}{{\bf{F}}})$ is a division algebra. The most famous example of a non-split quaternion algebra is Hamilton's quaternions  ${\mathbb{H}}=(\frac{-1,-1} {\mathbb{R}})$.
  \item $(\frac{a,b}{{\bf{F}}})$ is isomorphic to
 the algebra of all $2\times 2$ matrices with entries from ${\bf{F}}$. In this case we say
that the ${\bf{F}}$-algebra is split.
\end{enumerate}
In contrast to a quaternion division algebra,
a split quaternion algebra contains zero-divisors, nilpotent
elements and nontrivial idempotents.

One of the most famous split quaternion algebras is  the split quaternions of James Cockle (\cite{coc})  ${\bf{H}}_{\bf{S}}(\frac{-1,1}{{\rm {\mathbb{R}}}})$, which can be represented as
\[{\bf{H}}_{\bf{S}}=\{q=q_{0}+ q_{1}i+q_{2}j+q_{3}k:  \{q_{0}, q_{1}, q_{2}, q_{3}\}\in {\mathbb{R}}\}.\]
${\bf{H}}_{\bf{S}}$ is an associative,
non-commutative, non-division ring with four basis elements $\{1, i, j, k\}$ satisfying
the equalities
\begin{gather*}
  i^{2}=-1, \,
  j^{2}=k^{2}=1, \\
  ij=-ji=k, \,jk=-kj=-i,\,ik=-ki=-j.
\end{gather*}
The split quaternions of James Cockle are also named coquaternions.
In this paper we shall consider coquaternions and denote their  ${\bf{H}}$ to simplify.

Coquaternions is a recently developing topic. There are some studies related to geometric applications of split quaternions such as (\cite{ku}-\cite{br}). Particularly, the geometric and physical applications of coquaternions require solving coquaternionic equations and their systems. Therefore,
there are many studies on coquaternionic  equations. We mention only some recent papers. The method of rearrangements has been used to solve linear quaternionic equations involving $axb$ in \cite{sh}, new method of solving general linear coquaternionic equations with the terms of
the form $axb$ has been obtained in \cite{er}. The properties of coquaternion matrices  has been discussed in \cite{al}. Particularly, in \cite{al} the authors have defined the complex adjoint matrix of coquaternion matrices and   given the definition
of q-determinant of coquaternion matrices that is an usual determinant of the complex adjoint matrix.

Recently, in \cite{ky} the concept of  immanant (consequently, determinant and permanent) has been extended to a split quaternion algebra using methods of the theory of the row and column determinants. The theory of the row and column determinants was introduced
 in \cite{ky1, ky2} for matrices   over the quaternion non-split algebra.
This theory  over the quaternion skew field is being actively developed as by  the author (see, for ex.\cite{ky3}-\cite{ ky6}), and others (see, for ex. \cite{so1}-\cite{so3}).

In this paper  properties of the determinant of a Hermitian matrix over  ${\bf{H}}$ will been investigated, and determinantal representations   will been  given for the inverse of a Hermitian coquaternion matrix.
Firstly in Section
2, we shall give some properties of coquaternions, coquaternion matrices, and noncommutative determinants in Subsection 2.1,  and some basic concepts and
results from the theory of the row-column determinants of coquaternion matrices in Subsection 2.2. We shall consider the lemma  about expanding row and column determinants  by cofactors along corresponding rows and columns in this subsection as well.
In Section 3, properties of the determinant of a Hermitian coquaternion matrix will been investigated by using row-column determinants. In Section 4, determinantal representations for inverses of Hermitian coquaternion matrix will been  given and    Cramer' rules for left and right systems of linear equations will been obtained.
In Section 5, we shall get Cramer's rule for two-sided coquaternionic matrix equations ${\bf AXB}={\bf D}$, where ${\bf A}$, ${\bf B}$ are Hermitian. The main results will been illustrated by examples.

\section{Preliminaries}
\subsection{Coquaternions, coquaternion matrices and noncommutative determinants}

For any coquaternion $q=q_{0}+ q_{1}i+q_{2}j+q_{3}k\in {\bf{H}}$, by ${\rm Re}\,q:=q_{0}$ and ${\rm Im}\,q:=q_{1}i+q_{2}j+q_{3}k$, we define the real and imagine parts of $q$, respectively. The conjugate of a coquaternion $q$ is $\bar{q}=q_{0}- q_{1}i-q_{2}j-q_{3}k$, then the trace ${\bf t}(q)=2{\rm Re}\,q=2q_{0}$
and the norm form ${\bf n}(q)=q\bar{q}=q_{0}^{2}+ q_{1}^{2}-q_{2}^{2}-q_{3}^{2}$. The norm form of an coquaternion $q$ usually denote by $I_{q}:={\bf n}(q)$. The norm of a coquaternion by $\|q\|=\sqrt{|I_{q}|}$ are considered as well. If  $\|q\|=1$, then $q$ is called unit coquaternion. Notice that $p=\frac{q}{\|q\|}$ is a unit coquaternion for $q\in {\bf{H}}$ with $\|q\| \neq 0$ and 1, $i$, $j$ and $k$ are the
basis units. The inverse of
the coquaternion q is $q^{-1}=\frac{\bar{q}}{I_{q}}$, where $I_{q} \neq 0$. We indicate by $\mathcal{U}({\bf{H}})$ the set of all invertible
elements of ${\bf{H}}$ and $\mathcal{D}({\bf{H}})$ the set of all zero-divisors of ${\bf{H}}$.

Denote by ${\bf H}^{n \times m}$ a set
of $n \times m$ matrices with entries in ${\bf H}$ and by ${\rm M}\left( {n,{\bf H}}\right)$ a ring  of $n\times n$ coquaternionic matrices. This is a ring with a unit which is the usual identity matrix ${\bf I}_{n}$. By usual way, we define the transpose ${\bf A}^{T}=(a_{ji})\in{\bf H}^{n\times m}$, the conjugate $\overline{{\bf A}}=(\overline{{a}_{ij}})\in{\bf H}^{m\times n}$, the Hermitian adjoint matrix (the conjugate transpose)  $ {\rm {\bf A}}^{ *}=(\overline{{a}_{ji}})\in{\bf H}^{n\times m}$ of ${\rm {\bf A}}=(a_{ij})\in{\bf H}^{m\times n}$,  and the inverse $ {\rm {\bf A}}^{ -1}$ of ${\rm {\bf A}}=(a_{ij})\in{\bf H}^{n\times n}$.
For more properties of split quaternions the
reader is referred to \cite{er1}-\cite{cl}.

Definition of  determinant of matrices  with noncommutative entries (that are also defined as noncommutative determinants) is more associated with matrices over the skew field of Hamilton's quaternions ${\mathbb H}$. There are even three   approaches in its defining. The first approach  to defining
the determinant of a matrix in ${\rm M}\left( {n,{\mathbb H}} \right)$
is as follows \cite{as,co}.

 \begin{defn} \label{defin:axiom} Let a
functional ${\rm d}:{\rm M}\left( {n,{\mathbb H}} \right) \to \bf {\mathbb H}$
satisfy the following axioms.
\begin{itemize}
\item [] {\bf Axiom 1} ${\rm d}\left( {{\rm {\bf A}}} \right) = 0$ if and only if the matrix ${\rm {\bf A}}$ is non invertible.
\item [] {\bf Axiom 2} ${\rm d}\left( {{\rm {\bf A}} \cdot
{\rm {\bf B}}} \right) = {\rm d}\left( {{\rm {\bf A}}} \right)
\cdot {\rm d}\left( {{\rm {\bf B}}} \right)$ for $\forall {\rm
{\bf B}}\in {\rm M}\left( {n,{\mathbb H}} \right)$.
\item [] {\bf Axiom 3} If the matrix ${\rm {\bf A}}'$ is obtained
from ${\rm {\bf A}}$ by adding a left-multiple of a row to another
row or a right-multiple of a column to another column, then ${\rm
d}\left( {{\rm {\bf A}}}\right)'={\rm d}\left( {{\rm {\bf
A}}}\right)$.
\end{itemize}
Then  the functional ${\rm d}$ is called the determinant
of  ${\rm {\bf A}}\in {\rm M}\left( {n,{\mathbb H}}
\right)$.\end{defn}  But in \cite{as}, it is proved that if a determinant functional satisfies
Axioms 1, 2,
 3, then its value is real. The famous examples of such determinant are the
 determinants of Study and  Diedonn\'{e}.

 In another way of looking  a noncommutative determinant  is
defined as a rational function from  entries. In particular, in the theory of the Gelfand-Retah   quasideterminants  \cite{ge1,ge2}, an arbitrary $n\times n$
matrix over a skew field has been associated with an $n\times n$ matrix
whose entries are quasideterminants. The quasideterminant is not
an analog of the usual determinant but rather of a ratio of
the determinant of an $n\times n$-matrix to the determinant of an
$(n-1)\times (n-1)$-submatrix.

At last, at the third approach  a noncommutative determinant is
defined, by analogy to the usual determinant, as the alternating sum of $n!$ products of entries of a matrix but by specifying a certain ordering of coefficients in
each term.  Moore \cite{mo} was the first who
 achieved the fulfillment of the
main Axiom 1 by such definition of  a noncommutative determinant.
But it has been done  not for all square matrices over a skew field but
only Hermitian matrices.  Later, Dyson  \cite{dy} described the theory in more modern
terms. But until recently, the Moore determinant has not been extended  to arbitrary square matrices over ${\mathbb H}$.
 The full and natural extension of the definition of  Moore's determinant to arbitrary square matrices over ${\mathbb H}$ has been reached in the theory of column-row determinants.

Recently  in \cite{al} the $q$-determinant of coquaternionic matrices has been introduced by the follows.  Let ${\bf A}={\bf A}_1+{\bf A}_2j \in{\bf H}^{n\times n}$, where ${\bf A}_1$ and ${\bf A}_2$ are complex matrices. Then the complex adjoint matrix $\chi_{A} \in{\mathbb C}^{2n\times 2n}$ is defined as
 \[\chi_{A}:=\begin{pmatrix}
 {\bf A}_1 & {\bf A}_2\\
\overline{{\bf A}_2} & \overline{{\bf A}_1}
\end{pmatrix}\]
and the $q$-determinant of ${\bf A}$ is defined as the usual determinant of  $\chi_{A}$, that is $|{\bf A}|_{q}=|\chi_{A}|$. It has been shown that  properties of the $q$-determinant is   close to the usual determinant, especially, it satisfies Axioms 1, 2. Since the $q$-determinant of ${\bf A}\in{\bf H}^{n\times n}$ takes a value not in ${\bf H}$ but in ${\mathbb C}$ and the $q$-determinant can not be expanded  by cofactors along an
arbitrary row or column, then determinantal representations of the inverse ${\bf A}^{-1} $ by the $q$-determinant could not be obtained.

\subsection{Definitions and basic properties   of the column and  row determinants}
For  ${\rm {\bf A}}=(a_{ij})\in{\bf H}^{n\times n}$ we define $n$ row determinants as  follows.

\begin{defn}\label{def_rdet}\cite{ky}
The $i$-th row determinant of  ${\rm {\bf A}}=(a_{ij})\in{\bf H}^{n\times n}$ is defined as
\[
{\rm{rdet}}_{ i} {\rm {\bf A}} = {\sum\limits_{\sigma \in S_{n}}
{\left( { - 1} \right)^{n - r}{a_{i{\kern 1pt} i_{k_{1}}} }
{a_{i_{k_{1}}   i_{k_{1} + 1}}}   \ldots } } {a_{i_{k_{1} + l_{1}}
 i}}  \ldots  {a_{i_{k_{r}}  i_{k_{r} + 1}}}
\ldots  {a_{i_{k_{r} + l_{r}}  i_{k_{r}} }},
\]
where left-ordered cycle notation of the permutation $\sigma$ is written
as follows
\begin{equation*}
   \sigma = \left( {i\,i_{k_{1}}  i_{k_{1} + 1} \ldots i_{k_{1} +
l_{1}} } \right)\left( {i_{k_{2}}  i_{k_{2} + 1} \ldots i_{k_{2} +
l_{2}} } \right)\ldots \left( {i_{k_{r}}  i_{k_{r} + 1} \ldots
i_{k_{r} + l_{r}} } \right).
\end{equation*}
Here the index $i$ starts the first cycle from the left  and other
cycles satisfy the  conditions,
$
 i_{k_{2}}  < i_{k_{3}}  < \ldots < i_{k_{r}},
\quad i_{k_{t}}  < i_{k_{t} + s},
$
for all $t = \overline {2,r}$ and $s = \overline {1,l_{t}}$, (since  ${\rm sign}(\sigma)=\left( { - 1} \right)^{n - r}$).
\end{defn}
For ${\rm {\bf A}}=(a_{ij})\in{\bf H}^{n\times n}$ we  define $n$ column determinant as well.

\begin{defn}\label{def_cdet}\cite{ky}
The  $j$-th column determinant of  ${\rm {\bf A}}=(a_{ij})\in{\bf H}^{n\times n}$ is defined as
\[
{\rm{rdet}}_{ j} {\rm {\bf A}} = {\sum\limits_{\tau \in S_{n}}
\left( { - 1} \right)^{n - r}{{a_{j_{k_{r}} j_{k_{r} +
l_{r}} } \ldots a_{j_{k_{r} + 1}  j_{k_{r}} }  \ldots } }a_{j\,
j_{k_{1} + l_{1}} }  \ldots  a_{ j_{k_{1} + 1} j_{k_{1}}
}a_{j_{k_{1}} j}},
\]
where right-ordered cycle notation of the
permutation $\tau \in S_{n}$ is written as follows
\begin{equation*}
 \tau = \left( {j_{k_{r} + l_{r}}  \ldots j_{k_{r} + 1} j_{k_{r}}
} \right)\ldots \left( {j_{k_{2} + l_{2}}  \ldots j_{k_{2} + 1}
j_{k_{2}} } \right){\kern 1pt} \left( {j_{k_{1} + l_{1}}  \ldots
j_{k_{1} + 1} j_{k_{1} } j} \right).
\end{equation*}
Here the first cycle from the right begins with the index $j$ and
other cycles satisfy the following conditions,
$
 j_{k_{2}}  < j_{k_{3}}  < \ldots < j_{k_{r}},
\quad j_{k_{t}}  < j_{k_{t} + s},
$
for all $t = \overline {2,r}$ and $s = \overline {1,l_{t}}$
\end{defn}

In \cite{ky}
the basic properties of the column and row imanants of a
square matrix over ${\rm {\bf{H}}}$ has been consider. These properties can be evidently extend to column-row determinants.
\begin{prop}\label{theorem:zeros_col_row}(The first theorem about zero of an row-column determinant)
If one of the rows (columns) of  ${\rm {\bf A}}=(a_{ij})\in{\bf H}^{n\times n}$  consists of zeros only, then $
{\rm{rdet}}_{{i}}\, {\rm {\bf A}} = 0$ and ${\rm{cdet}} _{{i}}\,
{\rm {\bf A}} = 0$ for all ${i = \overline {1,n}}. $
\end{prop}

Denote by ${\bf{H}}a$ and $a{\bf{H}}$  left  and  right principal ideals of ${\bf{H}}$, respectively.
\begin{prop}\label{theorem:zeros_row2}(The second theorem about zero of an row determinant)
Let ${\rm {\bf A}}=(a_{ij})\in{\bf H}^{n\times n}$ and  $a_{ki}\in{\bf{H}}a_{i}$ and $a_{ij}\in \overline{a_{i}}{\bf{H}}$, where $n(a_{i})=0$ for $k,j = \overline {1,n}$ and  for all $i\neq k$. Let $a_{11}\in{\bf{H}}a_{1}$ and $a_{22}\in \overline{a_{1}}{\bf{H}}$ if $k=1$, and $a_{kk}\in{\bf{H}}a_{k}$ and $a_{11}\in \overline{a_{k}}{\bf{H}}$ if $k=i>1$, where $n(a_{k})=0$. Then ${\rm{rdet}}_{k} {\rm {\bf A}} =0$.
\end{prop}

\begin{prop}\label{theorem:zeros_col2}(The second theorem about zero of an column determinant)
Let ${\rm {\bf A}}=(a_{ij})\in{\bf H}^{n\times n}$ and  $a_{ik}\in a_{i}{\bf{H}}$ and $a_{ji}\in {\bf{H}}\overline{a_{i}}$, where $n(a_{i})=0$ for $k,j = \overline {1,n}$ and  for all $i\neq k$. Let $a_{11}\in a_{1}{\bf{H}}$ and $a_{22}\in {\bf{H}}\overline{a_{1}}$ if $k=1$, and $a_{kk} \in a_{k}{\bf{H}}$ and $a_{11}\in {\bf{H}}\overline{a_{k}}$ if $k=i>1$, where $n(a_{k})=0$. Then ${\rm{cdet}}_{k} {\rm {\bf A}} =0$.
\end{prop}

\begin{prop}\label{theorem:l_mult} If the $i$-th row of  ${\rm {\bf A}}=(a_{ij})\in{\bf H}^{n\times n}$ is left-multiplied by $b \in {\bf{H}} $, then $
{\rm{rdet}}_{{i}}\, {\rm {\bf A}}_{{i{\kern 1pt}.}} \left( {b
\cdot {\rm {\bf a}}_{{i{\kern 1pt}.}}} \right) = b \cdot
{\rm{rdet}}_{{i}}\, {\rm {\bf A}}$ for all ${i = \overline {1,n}}
.$
\end{prop}
\begin{prop}\label{theorem:c_mult} If the $j$-th column of
 ${\rm {\bf A}}=(a_{ij})\in{\bf H}^{n\times n}$ is right-multiplied by $b \in {\bf{H}}$, then $
{\rm{cdet}} _{{j}}\, {\rm {\bf A}}_{{.{\kern 1pt}j}} \left( {{\rm
{\bf a}}_{{.{\kern 1pt}j}} \cdot b} \right) = {\rm{cdet}} _{{j}}\,
{\rm {\bf A}}\cdot b$ for all ${j = \overline {1,n}}.$
\end{prop}
\begin{prop}\label{theorem:sum_matr_row} If for  ${\rm {\bf A}}=(a_{ij})\in{\bf H}^{n\times n}$ there exists  $t \in \{1,...,n\} $  such that $a_{tj} =
b_{j} + c_{j} $ for all $j = \overline {1,n}$, then for all $i =
\overline {1,n}$
\[
\begin{array}{l}
   {\rm{rdet}}_{{i}}\, {\rm {\bf A}} = {\rm{rdet}}_{{i}}\, {\rm {\bf
A}}_{{t{\kern 1pt}.}} \left( {{\rm {\bf b}}} \right) +
{\rm{rdet}}_{{i}}\, {\rm {\bf A}}_{{t{\kern 1pt}.}} \left( {{\rm
{\bf c}}} \right), \\
  {\rm{cdet}} _{{i}}\, {\rm {\bf A}} = {\rm{cdet}} _{{i}}\, {\rm
{\bf A}}_{{t{\kern 1pt}.}} \left( {{\rm {\bf b}}} \right) +
{\rm{cdet}}_{{i}}\, {\rm {\bf A}}_{{t{\kern 1pt}.}} \left( {{\rm
{\bf c}}} \right),
\end{array}
\]
where ${\rm {\bf b}}=(b_{1},\ldots, b_{n})\in{\bf H}^{1\times n}$, ${\rm {\bf
c}}=(c_{1},\ldots, c_{n})\in{\bf H}^{1\times n}$ are arbitrary row-vectors.
\end{prop}
\begin{prop}\label{theorem:sum_matr_col} If for ${\rm {\bf A}}=(a_{ij})\in{\bf H}^{n\times n}$ there exists  $t \in \{1,...,n\} $ such that
$a_{i\,t} = b_{i} + c_{i}$  for all  $i = \overline {1,n}$, then for all $j
= \overline {1,n}$
\[
\begin{array}{l}
  {\rm{rdet}}_{{j}}\, {\rm {\bf A}} = {\rm{rdet}}_{{j}}\, {\rm {\bf
A}}_{{.\,{\kern 1pt}t}} \left( {{\rm {\bf b}}} \right) +
{\rm{rdet}}_{{j}}\, {\rm {\bf A}}_{{.\,{\kern 1pt} t}} \left(
{{\rm
{\bf c}}} \right),\\
  {\rm{cdet}} _{{j}}\, {\rm {\bf A}} = {\rm{cdet}} _{{j}}\, {\rm
{\bf A}}_{{.\,{\kern 1pt}t}} \left( {{\rm {\bf b}}} \right) +
{\rm{cdet}} _{{j}} {\rm {\bf A}}_{{.\,{\kern 1pt}t}} \left( {{\rm
{\bf c}}} \right),
\end{array}
\]
where ${\rm {\bf b}}=(b_{1},\ldots, b_{n})^T\in{\bf H}^{n\times 1}$, ${\rm {\bf
c}}=(c_{1},\ldots, c_{n})^T\in{\bf H}^{n\times 1}$ are arbitrary column-vectors.
\end{prop}
\begin{prop}  If ${\rm {\bf A}}^{ *} $
is the Hermitian adjoint matrix (the conjugate transpose) of  ${\rm {\bf A}}=(a_{ij})\in{\bf H}^{n\times n}$, then $
{\rm{rdet}}_{{i}}\, {\rm {\bf A}}^{ *} =  \overline{{{\rm{cdet}}
_{{i}}\, {\rm {\bf A}}}}$ for all $i = \overline {1,n} $.
\end{prop}

The  following lemma enables  to expand
 ${\rm{rdet}}_{{i}}\, {\rm {\bf A}}$
 by cofactors along  the $i$-th row for all ${i = \overline {1,n}} $. Consequently, the calculation of the row
determinant of a $n\times n$ matrix is reduced to the calculation
of the row determinant of a lower dimension matrix.
\begin{defn} Let ${\rm {\bf A}} \in {\rm M}(n,{{\bf H}})$  and
$ {\rm{rdet}}_{{i}}\, {\rm {\bf A}} =  {\sum\limits_{j} {a_{ij}
} } R _{ i{\kern 1pt}j},$ for all $i = \overline
{1,n}$. Then $R _{ ij} $ is called
 the right  $ij$-th cofactor  of ${\rm {\bf A}}$.
\end{defn}

\begin{lem}\label{lemma:R_ij} Let $R_{ij}$ be the right $ij$-th
cofactor of ${\rm {\bf A}}\in {\rm M}\left( {n,{\bf H}}
\right)$, that is $ {\rm{rdet}}_{{i}}\, {\rm {\bf A}} =
{\sum\limits_{j = 1}^{n} {{a_{i{\kern 1pt} j} \cdot R_{i{\kern
1pt} j} } }} $ for all ${i = \overline {1,n}}$.  Then
\begin{equation}\label{eq:Rij}
 R_{i{\kern 1pt} j} ={\left\{\begin{array}{ccc}
                                - {\rm{rdet}}_{{k}}\, {\rm {\bf A}}_{{.{\kern 1pt} j}}^{{i i}} \left( {{\rm
{\bf a}}_{{. i}}}  \right), &{i \ne j}  & k=\left\{\begin{array}{ccc}
                                                     j, & if& i>j;\\
                                                     j-1,& if& i<j;
                                                   \end{array}\right.
\\
                               {\rm{rdet}} _{{k}}\, {\rm {\bf A}}^{{i i}}, &{i = j} &  k = \min {\left\{ {I_{n}}  \right.} \setminus {\left. {i}
\right\}}
                             \end{array}
 \right.}
\end{equation}
 where  ${\rm {\bf A}}_{. j}^{i i} \left(
{{\rm {\bf a}}_{. i}}  \right)$ is obtained
from ${\rm {\bf A}}$   by replacing the $j$-th column with the
$i$-th column, and then by deleting both the $i$-th row and column, $I_n=\{1,\ldots,n\}$.
\end{lem}
\emph{Proof.} At first  we  prove that $R_{i{\kern 1pt} i} =
{\rm{rdet}} _{{k}}\, {\rm {\bf A}}^{{i{\kern 1pt} i}}$, where $k =
\min {\left\{ {I_{n}}  \right.} \setminus {\left. {i}
\right\}} $.

 If $i = 1$, then ${\rm{rdet}} _{1}\, {\rm {\bf A}} =
{a_{11} \cdot R{}_{11} + a_{12} \cdot R{}_{12} + \ldots + a_{1n}
\cdot R{}_{1n}}$. Consider some monomial of ${\rm{rdet}} _{1}\,
{\rm {\bf A}}$ such that begin with $a_{1{\kern 1pt} 1} $ from the
left,
\begin{gather*}
   a_{11} \cdot R_{11} = {\sum\limits_{\tilde {\sigma}  \in S_{n}}
{\left( { - 1} \right)^{n - r}a_{11}  a_{2{\kern 1pt} i_{k_{2}} }
 \ldots } }a_{i_{k_{2} + l_{2}} {\kern 1pt} 2}  \ldots
 a_{i_{k_{r}} {\kern 1pt} i_{k_{r} + 1}}  \ldots
a_{i_{k_{r} + l_{r}}  {\kern 1pt} i_{k_{r}} }=\\ a_{11}  {\sum\limits_{\tilde {\sigma} _{1} \in S_{n-1} } {\left( { - 1} \right)^{n - 1 - \left( {r - 1}
\right)}a_{2{\kern 1pt}  i_{k_{2}} }  \ldots} } a_{i_{k_{2} +
l_{2}}   2} \ldots a_{i_{k_{r}}  i_{k_{r} + 1}} \ldots a_{i_{k_{r}
+ l_{r} } i_{k_{r}} },
\end{gather*}
where
\begin{gather*} \tilde
{\sigma}  = \left( {1} \right)\left( {2\,i_{k_{2}}  \ldots
i_{k_{2} + l_{2}} }  \right)\ldots \left( {i_{k_{r}}  i_{k_{r} +
1} \ldots i_{k_{r} + l_{r}} }  \right),\\ \tilde {\sigma} _{1} = \left( {2\,i_{k_{2}}  \ldots
i_{k_{2} + l_{2}} } \right)\ldots \left( {i_{k_{r}}  i_{k_{r} + 1}
\ldots i_{k_{r} + l_{r}} } \right).\end{gather*}
$S_{n-1}$ is the
symmetric group on the set $I_{n}\setminus 1$. The numbers of the
disjoint cycles and  coefficients of every monomial of $R_{11}$
decrease  by one. Since   elements of the second row start  these monomials   on the left and elements of the first row and  column do not belong to their, then
\begin{equation}\label{kyr3}
 R_{11} ={\sum\limits_{\tilde {\sigma} _{1} \in S_{n-1} } {\left( { - 1}
\right)^{n - 1 - \left( {r - 1} \right)}a_{2{\kern 1pt}  i_{k_{2}}
}  \ldots} } a_{i_{k_{2} + l_{2}}   2} \ldots a_{i_{k_{r} + l_{r}
}  i_{k_{r}} }  = {\rm{rdet}} _{2} {\rm {\bf A}}^{11}.
\end{equation}

If now $i \ne 1$, then
\begin{equation}\label{kyr4}
 {\rm{rdet}}_{{i}}\, {\rm {\bf A}} ={ a_{i1} \cdot R_{i1} +
a_{i2} \cdot R{}_{i2} + \ldots + a_{i{\kern 1pt} n} \cdot
R{}_{i{\kern 1pt} n}}
\end{equation}
Consider some monomial of ${\rm{rdet}}_{{i}}\, {\rm {\bf A}}$ such that begins with $a_{i{\kern 1pt} i} $ from the left,
\begin{gather*}
   a_{i{\kern 1pt} i} \cdot R_{i{\kern 1pt} i} =
{\sum\limits_{\mathord{\buildrel{\lower3pt\hbox{$\scriptscriptstyle\frown$}}\over
{\sigma} } \in \,S_{n}}  {\left( { - 1} \right)^{n - r}a_{i{\kern
1pt} i}  a_{1{\kern 1pt}  i_{k_{2}} }   \ldots } }a_{i_{k_{2} +
l_{2}}  1}  \ldots  a_{i_{k_{r}}  i_{k_{r} + 1}} \ldots
a_{i_{k_{r} + l_{r}}  i_{k_{r}} }=\\a_{i{\kern 1pt} i}
\cdot
{\sum\limits_{\mathord{\buildrel{\lower3pt\hbox{$\scriptscriptstyle\frown$}}\over
{\sigma} } _{1} \in
\mathord{\buildrel{\lower3pt\hbox{$\scriptscriptstyle\frown$}}\over{S}_{n-1}}}
{\left( { - 1} \right)^{n - 1 - \left( {r - 1} \right)}a_{1 {\kern
1pt} i_{k_{2} }}  \ldots} } a_{i_{k_{2} + l_{2}}   1}  \ldots
a_{i_{k_{r} + l_{r}}   i_{k_{r}} },
\end{gather*}
where \begin{gather*}
\mathord{\buildrel{\lower3pt\hbox{$\scriptscriptstyle\frown$}}\over
{\sigma }} = \left( {i} \right)\left( {1\,i_{k_{2}}  \ldots
i_{k_{2} + l_{2}} } \right)\ldots \left( {i_{k_{r}}  i_{k_{r} + 1}
\ldots i_{k_{r} + l_{r}} } \right),\\
\mathord{\buildrel{\lower3pt\hbox{$\scriptscriptstyle\frown$}}\over
{\sigma }} _{1} = \left( {1\,i_{k_{2}}  \ldots i_{k_{2} + l_{2}} }
\right)\ldots \left( {i_{k_{r}}  i_{k_{r} + 1} \ldots i_{k_{r} +
l_{r}} }  \right).\end{gather*}
$\mathord{\buildrel{\lower3pt\hbox{$\scriptscriptstyle\frown$}}\over{S}_{n-1}}$
is the symmetric group on $I_{n}\setminus
i$. The numbers of
disjoint cycles and the coefficients of every monomial of
$R_{i{\kern 1pt} i}$ again decrease  by  one. Each monomial of
$R_{i{\kern 1pt} i}$
  begins on the left with an entry of the first row.
Since   elements of the first row start  these monomials   on the left and elements of the $i$-th row and  column do not belong to their, then
\begin{equation}\label{kyr5}
R_{i{\kern 1pt} i}
={\sum\limits_{\mathord{\buildrel{\lower3pt\hbox{$\scriptscriptstyle\frown$}}\over
{\sigma} } _{1} \in
\mathord{\buildrel{\lower3pt\hbox{$\scriptscriptstyle\frown$}}\over{S}_{n-1}}}
{\left( { - 1} \right)^{n - 1 - \left( {r - 1} \right)}a_{1 {\kern
1pt} i_{k_{2} }}  \ldots} } a_{i_{k_{2} + l_{2}}   1} \ldots
a_{i_{k_{r} + l_{r}}   i_{k_{r}} }  = {\rm{rdet}} _{1} {\rm {\bf
A}}^{i{\kern 1pt} i}.
\end{equation}
By combining (\ref{kyr3}) and (\ref{kyr5}), we get $R_{i{\kern 1pt}
i} = {\rm{rdet}} _{{k}}\, {\rm {\bf A}}^{{i{\kern 1pt} i}}$, $k =
\min {\left\{ {I_{n} \setminus i} \right\}}$.

Now suppose that $i \ne j$. Consider  some monomial of ${\rm{rdet}}_{{i}}\, {\rm {\bf A}}$ in (\ref{kyr4})  such that begins with $a_{i{\kern 1pt} j} $ from the left,
\[
\begin{array}{c}
   a_{i{\kern 1pt} j} \cdot R_{i{\kern 1pt} j} = {\sum\limits_{\bar
{\sigma} \in\,  S_{n}}  {\left( { - 1} \right)^{n - r}a_{i{\kern
1pt} j}\,  a_{j{\kern 1pt}  i_{k_{1}} }  \ldots } } a_{i_{k_{1} +
l_{1}}   i}  \ldots  a_{i_{k_{r}}  i_{k_{r} + 1}} \ldots
a_{i_{k_{r} + l_{r}}  i_{k_{r}} }  =\\
    = - a_{i{\kern 1pt} j} \cdot {\sum\limits_{\bar {\sigma}  \in \,S_{n}}
{\left( { - 1} \right)^{n - r - 1}a_{j {\kern 1pt} i_{k_{1}} }
\ldots } } a_{i_{k_{1} + l_{1}}   i} \ldots a_{i_{k_{r}}  i_{k_{r}
+ 1}}   \ldots a_{i_{k_{r} + l_{r}}  i_{k_{r}} } ,\\
\end{array}
\]
where $ \bar {\sigma}  = \left( {i\,j\,\,i_{k_{1}}  \ldots
i_{k_{1} + l_{1}} } \right)\ldots \left( {i_{k_{r}}  i_{k_{r} + 1}
\ldots i_{k_{r} + l_{r}} } \right)$. Denote $\tilde {a}_{i_{k_{1}
+ l_{1}}  j} = a_{i_{k_{1} + l_{1}} i} $ for all $ i_{k_{1} +
l_{1}}  \in  {I_{n}}$. Then
\[
   a_{i{\kern 1pt} j} \cdot R_{i{\kern 1pt} j} =
 - a_{i{\kern 1pt} j} \cdot {\sum\limits_{\bar {\sigma} _{1} \in \,
 \mathord{\buildrel{\lower3pt\hbox{$\scriptscriptstyle\frown$}}\over{S}_{n-1}}}
 {\left( { - 1} \right)^{n - r -
1}a_{j{\kern 1pt} i_{k_{1}} }   \ldots } } \tilde {a}_{i_{k_{1} +
l_{1} }  j}    \ldots  a_{i_{k_{r} + l_{r}}   i_{k_{r}} },
\]
where $\bar {\sigma} _{1} = \left( {j\,i_{k_{1}}  \ldots i_{k_{1}
+ l_{1}} } \right)\ldots \left( {i_{k_{r}}  i_{k_{r} + 1} \ldots
i_{k_{r} + l_{r}} } \right)$. The permutation   $\bar {\sigma}
_{1} $  does not contain the index $i$ in each monomial of
$R_{i{\kern 1pt} j}$. This permutation satisfies the conditions of
Definition \ref{def_rdet} for  ${\rm rdet}_{j} {\rm {\bf A}}_{.{\kern
1pt} j}^{i{\kern 1pt} i} \left( {{\rm {\bf a}}_{. i}} \right)$. The matrix ${\rm {\bf A}}_{.
j}^{i{\kern 1pt} i} \left( {{\rm {\bf a}}_{. i}} \right) $  is obtained from ${\rm {\bf A}}$ by replacing
the $j$-th column with the column $i$, and then by deleting both
the $i$-th row and  column. That is,
\[
{\sum\limits_{\bar {\sigma} _{1} \in \,
 \mathord{\buildrel{\lower3pt\hbox{$\scriptscriptstyle\frown$}}\over{S}_{n-1}}}
  {\left( { - 1} \right)^{n - r - 1}a_{j{\kern 1pt}
 i_{k_{1}} }  \ldots} }\tilde {a}_{i_{k_{1} + l_{1}}
{\kern 1pt} j} \ldots a_{i_{k_{r} + l_{r}}  {\kern 1pt} i_{k_{r}}
}  = {\rm{rdet}} _{{j}}\, {\rm {\bf A}}_{{.j}}^{{i{\kern 1pt}i}}
\left( {{\rm {\bf a}}_{{.{\kern 1pt} i}}} \right)
\]
But ${\rm {\bf A}}_{.
j}^{i{\kern 1pt} i} \left( {{\rm {\bf a}}_{. i}} \right) $ is a quadratic matrix of order $n-1$. Therefore, more precisely on the set of indices of the matrix ${\rm {\bf A}}_{.{\kern 1pt}
j}^{i i} \left( {{\rm {\bf a}}_{. i}} \right)$  should be noted follows. If $i>j$,  then the index $j$ remains the same for ${\rm {\bf A}}_{.
j}^{i i} \left( {{\rm {\bf a}}_{. i}} \right)$ and
\begin{equation}\label{ind1}
R_{ij} = - {\rm rdet} _{j} {\rm
{\bf A}}_{. j}^{i{\kern 1pt} i} \left( {{\rm {\bf
a}}_{. i}} \right)
\end{equation}
But if  $i<j$,  then after deleting both
the $i$-th row and  column in $\bf A$ the  $j$-th row will be the  $j-1$-th row of ${\rm {\bf A}}_{.{\kern 1pt}
j}^{i{\kern 1pt} i} \left( {{\rm {\bf a}}_{.\,{\kern 1pt} {\kern
1pt} i}} \right)$. Therefore,
\begin{equation}\label{ind2}
R_{ij} = - {\rm rdet} _{j-1} {\rm
{\bf A}}_{. j}^{i i} \left( {{\rm {\bf
a}}_{. i}} \right)
\end{equation}
Combining (\ref{ind1}) and (\ref{ind2}), we finally obtain (\ref{eq:Rij}).
$\Box$
\begin{defn} Let ${\rm {\bf A}} \in {\rm M}(n,{{\bf H}})$  and
 $ {\rm cdet}_{j} {\rm {\bf A}} = {\sum\limits_{i}
{L _{ij} \, a_{ij}} },$ for all $j = \overline
{1,n} $. Then $L _{ij} $ is called the left
$ij$-th cofactor of ${\rm {\bf A}}$.
\end{defn}
\begin{lem}\label{lemma:L_ij} Let $L_{i{\kern 1pt} j} $ be
the left $ij$th cofactor of of a matrix ${\rm {\bf A}}\in {\rm
M}\left( {n,{\bf H}} \right)$, that is $ {\rm{cdet}}
_{{j}}\, {\rm {\bf A}} = {{\sum\limits_{i = 1}^{n} {L_{i{\kern
1pt} j} \cdot a_{i{\kern 1pt} j}} }}$ for all $j = \overline
{1,n}$. Then
\begin{equation}\label{eq:Lij}
L_{i{\kern 1pt} j} ={\left\{\begin{array}{lll}
                                - {\rm{cdet}}_{{k}}\, {\rm {\bf A}}_{{i{\kern 1pt}.}}^{{jj}} \left( {{\rm
{\bf a}}_{{j.}}}  \right), &{i \ne j}  & k=\left\{\begin{array}{lll}
                                                     i, & if& j>i;\\
                                                     i-1,& if& j<i;
                                                   \end{array}\right.
\\
                               {\rm{cdet}} _{{k}}\, {\rm {\bf A}}^{{i i}}, &{i = j} &  k = \min {\left\{ {J_{n}}  \right.} \setminus {\left. {j}
\right\}}
                             \end{array}
 \right.}
\end{equation}
where  ${\rm {\bf A}}_{i{\kern 1pt} .}^{jj} \left( {{\rm {\bf
a}}_{j{\kern 1pt} .} } \right)$ is obtained from ${\rm {\bf A}}$
 by replacing the $i$th row with the $j$th row, and then by
deleting both the $j$th row and  column, $J_n=\{1,\ldots,n\}$.
\end{lem}
\emph{Proof.} The proof is similar to the proof of Lemma \ref{lemma:R_ij}.
$\Box$

If ${\rm {\bf A}}^{ *}={\bf A}$, then ${\rm {\bf A}}\in{\bf H}^{n\times n}$ is called a Hermitian matrix. We finish this section by the following  theorem which is crucial for row-column determinants of a Hermitian matrix.
\begin{thm}\label{theorem:rdet_eq_cdet} If ${\rm {\bf A}}\in{\bf H}^{n\times n}$ is a Hermitian matrix, then
\[
{\rm{rdet}} _{1} {\rm {\bf A}} = \ldots = {\rm{rdet}} _{n} {\rm
{\bf A}} = {\rm{cdet}} _{1} {\rm {\bf A}} = \ldots = {\rm{cdet}}
_{n} {\rm {\bf A}}  \in {\mathbb{R}}.
\]
\end{thm}
By Theorem \ref{theorem:rdet_eq_cdet}, we have the following definition.
\begin{defn}\label{def_det__her}
Since all  column and  row determinants of a
Hermitian matrix over ${\bf{H}}$ are equal, we can
define the  determinant of a  Hermitian matrix ${\rm {\bf A}}\in{\bf H}^{n\times n}$. By definition, we
put for all $i = \overline {1,n}$,
\[
   \det {\rm {\bf A}}: = {\rm{rdet}}_{{i}}\,
{\rm {\bf A}} = {\rm{cdet}} _{{i}}\, {\rm {\bf A}}.
\]
\end{defn}
Evidently, if ${\bf A}\in{\bf{H}}^{2\times 2}$ is Hermitian and $a_{ij}\in\mathcal{D}({\bf{H}})$ for all $i,j=\overline {1,2}$, then $\det {\bf A}=0$. It would be expected in the general case, but the following example claims that it is not true.
\begin{ex}
Consider the Hermitian matrix \begin{equation}
\label{eq:A}{\rm {\bf A}}=\begin{pmatrix}
  0 & 1-k& 1-j \\
  1+k& 0  & 1+j\\
  1+j& 1-j  & 0
\end{pmatrix}.\end{equation}
 It can easily be checked that $a_{ij}\in\mathcal{D}({\bf{H}})$ for all $i,j=\overline {1,3}$. So,
 \begin{multline*}\det {\bf A}={\rm{rdet}}_{{1}}{\bf A}=\\ 0-0(1+j)(1-j)+(1-k)(1+j)(1+j)-(1-k)(1+k)0+\\(1-j)(1-j)(1+k)-(1-j)(1+j)0=4. \end{multline*}
\end{ex}

\section{ Properties of the column and  row  determinants of a Hermitian matrix}
\begin{thm} \label{repl_row} If the matrix ${\rm {\bf A}}_{j.} \left( {{\rm {\bf a}}_{i.}}
\right)$ is obtained from a Hermitian matrix ${\rm {\bf A}}\in
M\left( {n,{\bf H}} \right)$  by replacing  its $j$-th
row with the $i$-th row,  then for all $i,j=\overline{1,n}$ such
that $i\neq j$ we have
\begin{equation}
\label{eq:rdet0}{\rm{rdet}} _{{j}} {\rm {\bf A}}_{{j\,.}} \left( {{\rm {\bf
a}}_{{i\,.}}} \right) = 0.\end{equation}
\end{thm}
\emph{Proof.} We assume $n > 3$ for  ${\rm {\bf A}}\in {\rm
M}\left( {n,{\bf H}}\right)$. The case $n \le 3$ is
easily proved by a simple check. Consider some monomial $d$ of
${\rm{rdet}} _{{j}}\, {\rm {\bf A}}_{{j\,.}} \left( {{\rm {\bf
a}}_{{i\,.}}} \right)$. Suppose the index permutation of its
coefficients forms a direct product of $r$ disjoint cycles, and
denote $i = i_{s} $. Consider all
 possibilities of disposition of an entry of the $i_{s} $-th
row  in the monomial $d$.

(i) Suppose an entry of the $i_{s} $-th row is placed in $d$ such
that the index $i_{s} $ starts some disjoint cycle, i.e.:
\begin{equation}
\label{kyr9} d = ( - 1)^{n - r}a_{j{\kern 1pt} i_{1}}   \ldots
 a_{i_{k} j} \, u_{1} \ldots  u_{\rho}  \, a_{i_{s} i_{s +
1}}   \ldots  a_{i_{s + m} i_{s}}  \, v_{1} \ldots v_{p}
\end{equation}
\noindent Here we denote by $u_{\tau}  $ and $v_{t} $  products of
coefficients whose indices  form some disjoint cycles for all $
\tau = \overline {1,\rho}$ and $t = \overline {1,p}$ such that
$\rho + p = r - 2$ or there are no such products. For $d$ there
are the following three monomials  of ${\rm{rdet}} _{{j}}\, {\rm
{\bf A}}_{{j\,.}} \left( {{\rm {\bf a}}_{i\,.}} \right)$.
\begin{gather*}  d_{1} = ( - 1)^{n - r}a_{j{\kern 1pt} i_{1}}   \ldots a_{i_{k} j}
\, u_{1}  \ldots  u_{\rho}  \, a_{i_{s} i_{s + m}} \ldots  a_{i_{s
+ 1} i_{s}}\, v_{1} \ldots v_{p},\\
   d_{2} = ( - 1)^{n - r + 1}a_{j{\kern 1pt} i_{s + 1}} \ldots
 a_{i_{s + m} i_{s}}  \, a_{i_{s} i_{1}}   \ldots
a_{i_{k} j} \, u_{1}  \ldots u_{\rho}\, v_{1}
 \ldots v_{p},\\
   d_{3} = ( - 1)^{n - r + 1}a_{j{\kern 1pt} i_{s + m}}   \ldots
 a_{i_{s + 1} i_{s}}  a_{i_{s}{\kern 1pt} i_{1}} \ldots
 a_{i_{k} j} \, u_{1} \ldots  u_{\rho}  \, v_{1} \ldots
v_{p}.
\end{gather*}
Suppose $a_{j i_{1}}  \ldots  a_{i_{k} j} = x$ and $a_{i_{s} i_{s
+ 1}} \ldots a_{i_{s + m} i_{s}}  = y$, then $\overline {y} =
a_{i_{s} i_{s + m}}  \ldots a_{i_{s + 1} \,i_{s}} $. Taking into
account  $a_{j{\kern 1pt} i_{1}}  = a_{i_{s} i_{1}}  $,
$a_{j{\kern 1pt}i_{s - 1}}  = a_{i_{s} i_{s - 1}}$ and $
a_{j{\kern 1pt}i_{s + 1}} = a_{i_{s} i_{s + 1}}  $, we consider
the sum of these monomials.
\begin{multline}\label{eq:sum1}
d + d_{1} + d_{2} + d_{3} = ( - 1)^{n - r}(x  u_{1} \ldots
u_{\rho} \, y + x  u_{1}  \ldots u_{\rho} \overline {y}
 - y  x  u_{1} \ldots u_{\rho}-\\ \overline {y} \cdot x
u_{1} \ldots  u_{\rho})  v_{1}  \ldots v_{p}
  =( - 1)^{n - r} ( xu_{1} \ldots u_{\rho}
{\bf t}(y)-{\bf t}(y)xu_{1} \ldots u_{\rho}) v_{1}  \ldots  v_{p} =0.
\end{multline}
Thus among the monomials of ${\rm{rdet}} _{{j}}\, {\rm {\bf
A}}_{{j\,.}} \left( {{\rm {\bf a}}_{{i\,.}}}  \right)$  we find
 three monomials  for $d$ such  that the sum of these  monomials and
 $d$ is equal to zero.

If in (\ref{kyr9}) $m = 0$ or $m = 1$, we accordingly get such monomials,
\begin{gather*}
   \tilde {d} = ( - 1)^{n - r}a_{j{\kern 1pt} i_{1}}  \ldots
\cdot a_{i_{k} j} \, u_{1} \ldots  u_{\rho}  \, a_{i_{s} i_{s}} \,
v_{1} \ldots  v_{p},\\
    \mathord{\buildrel{\lower3pt\hbox{$\scriptscriptstyle\frown$}}\over {d}} = (
- 1)^{n - r}a_{j i_{1}}   \ldots  a_{i_{k} j} \, u_{1}
\ldots  u_{\rho}  \, a_{i_{s} i_{s + 1}} \, a_{i_{s + 1} i_{s}} \,
v_{1}  \ldots  v_{p}.
\end{gather*}
For them, there are the following  monomials, respectively,
\begin{gather*}
     \tilde {d}_{1} = ( - 1)^{n - r + 1}a_{j{\kern 1pt}i_{s}}  \, a_{i_{s}
{\kern 1pt} i_{1}}  \, \ldots  a_{i_{k} j} \, u_{1}
 \ldots  u_{\rho}  \, v_{1}  \ldots v_{p},\\
   \mathord{\buildrel{\lower3pt\hbox{$\scriptscriptstyle\frown$}}\over {d}}
_{1} = ( - 1)^{n - r + 1}a_{j{\kern 1pt}i_{s + 1}} \, a_{i_{s + 1}
{\kern 1pt} i_{s}}  \,a_{i_{s} {\kern 1pt} i_{1}} \ldots a_{i_{k}
j}\, u_{1} \ldots  u_{\rho}\,
 v_{1}  \ldots  v_{p}.
\end{gather*}
Taking into account $a_{j{\kern 1pt} i_{1}}  = a_{i_{s} i_{1}},$
$a_{j{\kern 1pt} i_{s}}  = a_{i_{s} i_{s}}\in {\rm {\mathbb{R}}}$,
$a_{j{\kern 1pt} i_{s + 1}} = a_{i_{s} i_{s + 1}} $, and
 $ {a_{i_{s} i_{s + 1}} a_{i_{s + 1} i_{s}}}= {\bf n}({a_{i_{s} i_{s + 1}}})
\in {\rm {\mathbb{R}}}$, we obtain $\tilde {d} + \tilde {d}_{1} = 0,$
\,
$\mathord{\buildrel{\lower3pt\hbox{$\scriptscriptstyle\frown$}}\over
{d}} +
\mathord{\buildrel{\lower3pt\hbox{$\scriptscriptstyle\frown$}}\over
{d}} _{1} = 0$. Hence,  the sums of corresponding two monomials of
${\rm{rdet}} _{{j}}\, {\rm {\bf A}}_{{j\,.}} ( {{\rm {\bf
a}}_{{i\,.}}})$ are equal to zero in these both cases.

 ii) Now suppose that the index $i_{s} $  is placed in  another disjoint
  cycle  than  $j$ and does not start this cycle,
\[
  \mathord{\buildrel{\lower3pt\hbox{$\scriptscriptstyle\smile$}}\over
{d}} = ( - 1)^{n - r}a_{j i_{1}}  \ldots a_{i_{k} j}\, u_{1}
\ldots u_{\rho }\, a_{i_{q} i_{q + 1}}  \ldots  a_{i_{s - 1}
i_{s}}
   a_{i_{s} i_{s + 1} } \ldots a_{i_{q - 1} i_{q}}  v_{1}
\ldots v_{p}.
\]
\noindent Here we denote by $u_{\tau}   $ and $v_{t} $  products
of coefficients whose indices   form some disjoint cycles for all
 $\tau = \overline {1,\rho}$ and $t = \overline
{1,p}$ such that $\rho + p = r - 2$ or there are no such products.
Now for $d$ there are the following three monomials   of
${\rm{rdet}} _{{j}}\, {\rm {\bf A}}_{{j\,.}} ( \rm {\bf
a}_{{i\,.}} )$,
\begin{gather*} \mathord{\buildrel{\lower3pt\hbox{$\scriptscriptstyle\smile$}}\over {d}}
_{1} = ( - 1)^{n - r}a_{j{\kern 1pt} i_{1}}  \ldots a_{i_{k} j}\,
u_{1} \ldots u_{\rho} \, a_{i_{q} i_{q - 1}}  \ldots a_{i_{s + 1}
i_{s}} a_{i_{s} i_{s - 1}}  \ldots \cdot a_{i_{q + 1} i_{q}} v_{1}
\ldots v_{p},\\
     \mathord{\buildrel{\lower3pt\hbox{$\scriptscriptstyle\smile$}}\over {d}}
_{2} = ( - 1)^{n - r + 1}a_{j{\kern 1pt} i_{s - 1}}  \ldots
a_{i_{q + 1} i_{q}}\,  a_{i_{q} i_{q - 1}}  \ldots a_{i_{s + 1}
i_{s}} a_{i_{s} i_{1}} \ldots a_{i_{k} j}\, u_{1} \ldots u_{\rho}
v_{1} \ldots v_{p},\\
\mathord{\buildrel{\lower3pt\hbox{$\scriptscriptstyle\smile$}}\over
{d}} _{3} = ( - 1)^{n - r + 1}a_{j{\kern 1pt} i_{s + 1}}  \ldots
a_{i_{q - 1} i_{q}} \, a_{i_{q} i_{q + 1}}  \ldots a_{i_{s - 1}
i_{s}} a_{i_{s} i_{1}} \ldots a_{i_{k} j} \, u_{1} \ldots u_{\rho}
v_{1} \ldots v_{p}.
\end{gather*}
Assume that \begin{gather*}a_{i_{s} i_{s + 1}}   \ldots \ a_{i_{q - 1} i_{q}}=
\varphi,\,a_{i_{q} i_{q + 1}}  \ldots  a_{i_{s - 1} i_{s}}=
\phi,\,a_{j\, i_{1}} \ldots a_{i_{k} j} =x,\\ a_{i_{q} i_{q +
1}} \ldots  a_{i_{s - 1} i_{s}} a_{i_{s} i_{s + 1} } \ldots
a_{i_{q - 1} i_{q}}  = y,\, a_{i_{s} i_{s + 1}} \ldots a_{i_{q
- 1} i_{q}} a_{i_{q} i_{q + 1}} \ldots a_{i_{s - 1} i_{s}}
=y_{1}.\end{gather*}
 Then we
obtain $y = \phi  \varphi $,\,  $y_{1} = \varphi \phi $,
$\overline {y} = a_{i_{q} i_{q - 1}}  \ldots a_{i_{s + 1} i_{s}}
a_{i_{s} i_{s - 1}}  \ldots a_{i_{q + 1} i_{q}}  $,\,\,and
$\overline {y_{1}} = a_{i_{s} i_{s - 1}} \ldots a_{i_{q + 1}
i_{q}} a_{i_{q} i_{q - 1}} \ldots a_{i_{s + 1} i_{s}}$. Accounting
for  $a_{j\, i_{1} } = a_{i_{s} i_{1}} $, $a_{ji_{s - 1}} =
a_{i_{s} i_{s - 1}}$, $a_{ji_{s + 1}} = a_{i_{s} i_{s + 1}}  $, we
have
\begin{multline}\label{eq:sum2}
\mathord{\buildrel{\lower3pt\hbox{$\scriptscriptstyle\smile$}}\over
{d}} +
\mathord{\buildrel{\lower3pt\hbox{$\scriptscriptstyle\smile$}}\over
{d}} _{1} +
\mathord{\buildrel{\lower3pt\hbox{$\scriptscriptstyle\smile$}}\over
{d}} _{2} +
\mathord{\buildrel{\lower3pt\hbox{$\scriptscriptstyle\smile$}}\over
{d}} _{3} =\\
 =( - 1)^{n - r}(xu_{1}  \ldots  u_{\rho} y+
  xu_{1}  \ldots  u_{\rho}  \overline {y}
- y_{1} xu_{1}  \ldots  u_{\rho}  - \overline {y_{1}}\, x
  u_{1}  \ldots  u_{\rho} )\times\\ \times v_{1}
\ldots  v_{p}= ( - 1)^{n - r}( xu_{1}  \ldots  u_{\rho}  {\bf t}(y) -
{\bf t}(y_{1} ) x
 u_{1} \ldots
u_{\rho})v_{1}  \ldots  v_{p}  = \\
 =  ( - 1)^{n - r}({\bf t}(\phi \cdot \varphi ) - {\bf t}(\varphi \cdot
\phi ))xu_{1}  \ldots u_{\rho} v_{1}  \ldots v_{p}.
\end{multline}
Since by the rearrangement property of the trace, ${\bf t}(\phi \cdot
\varphi ) = t(\varphi \cdot \phi )$, then we obtain
$\mathord{\buildrel{\lower3pt\hbox{$\scriptscriptstyle\smile$}}\over
{d}} +
\mathord{\buildrel{\lower3pt\hbox{$\scriptscriptstyle\smile$}}\over
{d}} _{1} +
\mathord{\buildrel{\lower3pt\hbox{$\scriptscriptstyle\smile$}}\over
{d}} _{2} +
\mathord{\buildrel{\lower3pt\hbox{$\scriptscriptstyle\smile$}}\over
{d}} _{3} = 0$.

(iii) If the indices  $i_{s} $ and $j$ are placed  in the same
cycle, then we have the following monomials: $d_{1} ,\;\tilde
{d}_{1}
,\;\mathord{\buildrel{\lower3pt\hbox{$\scriptscriptstyle\frown$}}\over
{d}} _{1}$ or
$\mathord{\buildrel{\lower3pt\hbox{$\scriptscriptstyle\smile$}}\over
{d}} _{1} $. As shown above,   for each of them  there are another
one or three  monomials of ${\rm{rdet}} _{{j}}\, {\rm {\bf
A}}_{{j\,.}} ( {\rm {\bf a}}_{{i{\kern 1pt}.}} )$ such that the
sums of these two  or  four corresponding monomials are equal to
zero.

We have considered all possible kinds of disposition of an entry
of the $i$-th row as a factor of some monomial $d$ of
${\rm{rdet}} _{{j}}\, {\rm {\bf A}}_{{j\,.}} ( {\rm {\bf
a}}_{{i{\kern 1pt}.}})$. For $d$, in each case there exist  one or three
corresponding monomials   such that accordingly the sum of the two or
four monomials is equal to zero. Thus, we have (\ref{eq:rdet0}).

We note that if one of factors of $d$ is zero, then evidently $d$, $d_{1}$, $d_{2}$, $d_{3}$ are equal $0$.
If two adjacent factors of $d$ are adjoint zero divisors (i.e. their product equals zero),
 then the sums (\ref{eq:sum1}) or (\ref{eq:sum2}) contain these adjacent zero divisors as well. Hence, the sums will be equal zero by the same cause.
$\Box$

The following theorem can be proved similarly.
\begin{thm} \label{repl_col}If the matrix ${\rm {\bf A}}_{.{\kern 1pt} i}
\left( {{\rm {\bf a}}_{.j}}  \right)$ is obtained from a Hermitian
matrix ${\rm {\bf A}}\in {\rm M}\left( {n,{\bf H}}
\right)$ by replacing of its $i$-th column with the $j$-th column,
then ${\rm{cdet}} _{{i}}\, {\rm {\bf A}}_{{.\,i}} ( {\rm {\bf
a}}_{{.j}}) = 0$ for all $i,j=\overline{1,n}$ such that $i\neq j.$
\end{thm}
\begin{cor}\label{cor:same_rows}
If a Hermitian matrix ${\rm {\bf A}}\in {\rm M}\left( {n,{\bf H}} \right)$ consists two same rows (columns), then
$\det{\rm {\bf A}}=0$.
\end{cor}
\emph{Proof.} Suppose the $i$-th row of ${\rm {\bf A}}$
coincides with the $j$-th row, i.e. $a_{ik} = a_{jk} $ for all $k
\in I_{n} $ and ${\left\{ {i,j} \right\}} \in I_{n} $ such that $i
\ne j$. Then $\overline {a_{ik}}  = \overline {a_{jk}}  $ for all
$k \in I_{n}$. Since  ${\rm {\bf A}}$ is Hermitian, then $a_{ki} = a_{kj}$
 for all $ k \in I_{n}$, where ${\left\{
{i,j} \right\}} \in I_{n} $ and $i \ne j$. It means that ${\rm {\bf A}}$ has two same
corresponding columns as well. The matrix ${\rm {\bf A}}$ may be
represented as ${\rm {\bf A}}_{i.} \left( {{\rm {\bf a}}_{j.} }
\right)$, where  ${\rm {\bf A}}_{i.} \left( {{\rm {\bf
a}}_{j.} } \right)$ is obtained from ${\rm {\bf A}}$ by replacing
the $i$-th row with the $j$-th row. By Theorem \ref{repl_row},
we have,
$
\det {\rm {\bf A}} = {\rm{rdet}} _{i} {\rm {\bf A}} = {\rm{rdet}}
_{i} {\rm {\bf A}}_{i.} \left( {{\rm {\bf a}}_{j.}}  \right) =
0.
$
$\Box$

We are needed by the following lemmas.
\begin{lem}\cite{ky1} \label{lemma:sum_factors}  Let $T_{n} $ be the sum
of all  possible products  of  the $n$ factors, each of which are
either $h_{i} \in {\bf H}$ or $ \overline {h_{i}} $ for all $ i =
\overline {1,n}$, by specifying the ordering in the terms, $ T_{n} =
h_{1} \cdot h_{2} \cdot \ldots \cdot h_{n} + \overline {h_{1}} \cdot
h_{2} \cdot \ldots \cdot h_{n} + \ldots + \overline {h_{1}}  \cdot
\overline {h_{2}}  \cdot \ldots \cdot \overline {h_{n}}. $ Then
$T_{n}$ consists of the  $2^{n}$ terms and $T_{n} = {\bf t}\left(
{h_{1}} \right)\;{\bf t}\left( {h_{2}} \right)\;\ldots \;{\bf t}\left( {h_{n}} \right).$
\end{lem}
\begin{lem}\label{lemma: rdetA_mult_col_b}
If the matrix ${\rm {\bf A}}_{.{\kern 1pt}i} \left( {{\rm {\bf
a}}_{.{\kern 1pt} i}}\cdot b \right)$ is obtained from a Hermitian
matrix ${\rm {\bf A}}\in {\rm M}\left( {n,{\bf H}}
\right)$ by right-multiplying of its $i$-th column  by $b \in {\bf H}$,  then for all $i=\overline{1,n}$ we have
 $ {\rm{rdet}}_{i} \,{\rm {\bf A}}_{.{\kern 1pt}i}
\left({{\rm {\bf a}}_{.{\kern 1pt} i}}\cdot b\right)= \det {\rm {\bf
A}}\cdot b.$
\end{lem}
\emph{Proof.} Consider some monomial $d$ of ${\rm {\bf
A}}_{.{\kern 1pt}i} \left( {{\rm {\bf a}}_{.{\kern 1pt} i}}\cdot b
\right)$ for $i=\overline {1\mbox{,}n} $. Denote $i_{k_1 }
\mbox{:}=i$. Then,
\[
\begin{array}{c}
  d=(-1)^{n-r}a_{i_{k_1 } i_{k_1 +1} } \ldots {\kern 1pt}{\kern
1pt}a_{i_{k_1 +l_1 } i_{k_1 } }b{\kern 1pt}{\kern 1pt}{\kern
1pt}a_{i_{k_2 } i_{k_2 +1} } \ldots a_{i_{k_2 +l_2 } i_{k_2 } }
\ldots \times \\
  \times a_{i_{k_r } i_{k_r +1} } \ldots a_{i_{k_r +l_r
} i_{k_r } } =(-1)^{n-r}h_1 \cdot b\cdot h_2 \cdot \ldots \cdot h_r
,
\end{array}
\]
where $h_s =a_{i_{k_s } i_{k_s +1} } \cdot \ldots \cdot a_{i_{k_s
+l_s } i_{k_s } } $ for all ${\ s=\overline {1,r} } $. If $l_s =1$,
then $h_s =a_{i_{k_s } i_{k_s +1} } \cdot a_{i_{k_s +1{\kern 1pt}}
i_{k_s } } ={\bf n}(a_{i_{k_s } i_{k_s +1} } )\in {\rm {\mathbb{R}}}$, and
if $l_s =0$, then $h_s =a_{i_{k_s } i_{k_s } } \in {\mathbb{R}}$.
Suppose there exists  such $s$ that $l_s \ge 2$. By Definition
\ref{def_rdet}, the index permutation $\sigma $ of $d$ forms a
direct products of disjoint cycles and its cycle notation is
left-ordered. Denote by $\sigma _s \left( {i_{k_s } }
\right)\mbox{:}=\left( {i_{k_s } i_{k_s +1} \ldots i_{k_s +l_s } }
\right)$ a cycle which corresponds to a factor $h_s $. Then $\sigma
_s^{-1} \left( {i_{k_s } } \right)\mbox{:}=\left( {i_{k_s } i_{k_s
+l_s } i_{k_s +1} \ldots i_{k_s +1} } \right)$ is the cycle which is
inverse to $\sigma _s \left( {i_{k_s } } \right)$ and corresponds to
the factor $\overline {h_s } $. There exist $2^{p-1}$ monomials of
${\rm {\bf A}}_{.{\kern 1pt}i} \left( {{\rm {\bf a}}_{.{\kern 1pt}
i}}\cdot b \right)$ such that their indices permutations form the
direct products of the disjoint cycles $\sigma _s \left( {i_{k_s } }
\right)$ or $\sigma _s^{-1} \left( {i_{k_s } } \right)$ for all
 $\left( {\ s=\overline {1,r} }
\right)$ and keeping their ordering from $1$ to $r$, where $p=r-\rho
$, and $\rho $ is the number of the cycles of the first and second
orders. Then by Lemma \ref{lemma:sum_factors} for the sum $C_1 $ of
these monomials and $d$ we obtain,
\[
C=(-1)^{n-r}b\cdot \alpha t(h_{\nu _1 } )\;\ldots \;t(h_{\nu _p }
),
\]
where $\alpha \in {\mathbb R}$ is a product of  factors whose
indices form  cycles of the first and second orders. Since ${\bf t}(h_{\nu _k } )\in {\mathbb R}$ for all $\nu _k \in \{1,\ldots
,r\}$ and $k=\overline {1,p} $, then $b$ commutes with ${\bf t}(h_{\nu _k } )\in {\mathbb R}$ for all $\nu _k \in \{1,\ldots
,r\}$ and $k=\overline {1,p} $. Therefore, ${\rm rdet} _i {\rm {\bf
A}}_{.\,i} \left( {{\rm {\bf a}}_{.\,i} \cdot b} \right)={\rm rdet}
_i {\rm {\bf A}}\cdot b=b\cdot \det {\rm {\bf A}}$.
$\Box$
\begin{lem}\label{lemma:
cdetA_mult_row_b} If ${\rm {\bf A}}_{i\,.} \left(b\cdot  { {\rm {\bf
a}}_{i\,.} }\right)$ is obtained from Hermitian ${\rm {\bf A}}\in
{\rm M}\left( {n,{\bf H}} \right)$ by left-multiplying of
its $i$-th row  by $b \in {\bf H}$, then for all
$i=\overline{1,n}$ we have
\[
 {\rm cdet}_i {\rm {\bf A}}_{i\,.} \left( {b\cdot {\rm {\bf a}}_{i\,.} }
\right)=b\cdot \det {\rm {\bf A}}\]
\end{lem}
The proof is similar to the proof of Lemma \ref{lemma:
rdetA_mult_col_b}.

By Theorems \ref{repl_row}, \ref{repl_col}, Lemmas  \ref{lemma:
rdetA_mult_col_b} and \ref{lemma: cdetA_mult_row_b}, and basic
properties of the row and column determinants, we have the following
theorems.
\begin{thm}\label{theorem:row_combin} If the $i$-th row of
a Hermitian matrix ${\rm {\bf A}}\in {\rm M}\left( {n,{\bf H}} \right)$ is replaced with a left linear combination
of its other rows, i.e. ${\rm {\bf a}}_{i.} = c_{1} {\rm {\bf
a}}_{i_{1} .} + \ldots + c_{k}  {\rm {\bf a}}_{i_{k} .}$, where $
c_{l} \in {{\bf H}}$ for all $ l = \overline {1,k}$ and
$\{i,i_{l}\}\subset I_{n} $, then
\[
 {\rm{rdet}}_{i}\, {\rm {\bf A}}_{i \, .} \left(
{c_{1} {\rm {\bf a}}_{i_{1} .} + \ldots + c_{k} {\rm {\bf
a}}_{i_{k} .}}  \right) = {\rm{cdet}} _{i}\, {\rm {\bf A}}_{i\, .}
\left( {c_{1}
 {\rm {\bf a}}_{i_{1} .} + \ldots + c_{k} {\rm {\bf
a}}_{i_{k} .}}  \right) = 0.
\]
\end{thm}
\begin{thm}\label{theorem:colum_combin} If the $j$-th column of
 a Hermitian matrix ${\rm {\bf A}}\in
{\rm M}\left( {n,{\bf H}} \right)$   is replaced with a
right linear combination of its other columns, i.e. ${\rm {\bf
a}}_{.j} = {\rm {\bf a}}_{.j_{1}}   c_{1} + \ldots + {\rm {\bf
a}}_{.j_{k}} c_{k} $, where $c_{l} \in{{\bf H}}$ for
all $l = \overline {1,k}$ and $\{j,j_{l}\}\subset J_{n}$, then
 \[{\rm{cdet}} _{j}\, {\rm {\bf A}}_{.j}
\left( {{\rm {\bf a}}_{.j_{1}} c_{1} + \ldots + {\rm {\bf
a}}_{.j_{k}}c_{k}} \right) ={\rm{rdet}} _{j} \,{\rm {\bf A}}_{.j}
\left( {{\rm {\bf a}}_{.j_{1}}  c_{1} + \ldots + {\rm {\bf
a}}_{.j_{k}}  c_{k}} \right) = 0.
\]
\end{thm}
 \begin{defn} Let ${\bf a}_{i.}\in{\bf H}^{n\times 1}$ for all $i=\overline{1,m}$.
Row-vectors $
 {\rm {\bf a}}_{1.},\ldots ,
 {\rm {\bf a}}_{m.} $
are left linearly dependent, if there exist scalars
$\{b_1,\ldots,b_m \}\subset {\bf H}$ (which are not all zero) such
that $ b_1 \cdot {\rm {\bf a}}_{1.} +\ldots +b_m
\cdot {\rm {\bf a}}_{m.} ={\rm {\bf 0}}$,
where  ${\rm {\bf 0}}$ is the zero row vector. If no such scalars exist, then the vectors are said to be left-linearly independent.
\end{defn}
  \begin{defn}Let ${\bf a}_{.j}\in{\bf H}^{1\times n}$ for all $j=\overline{1,m}$.
Column-vectors
$ {\rm {\bf a}}_{.1}, \ldots ,
 {\rm {\bf a}}_{.\,m} $
are   right linearly dependent, if there exist scalars
$\{c_1,\ldots,c_m \}\subset {\bf H}$ (which are not all zero) such
that ${\rm {\bf a}}_{.1 } \cdot  c_1 +\ldots +{\rm
{\bf a}}_{.m } \cdot c_m ={\rm {\bf 0}}$,
where  ${\rm {\bf 0}}$ is the zero column-vector. If no such scalars exist, then the column-vectors are said to be right-linearly independent.
\end{defn}
By Lemma \ref{cor:same_rows},  the  evident corollary of Theorems  \ref{theorem:row_combin} and \ref{theorem:colum_combin} follows.
\begin{cor}\label{cor:colum_combin} If the $i$-th row of
  Hermitian  ${\rm {\bf A}}\in
{\rm M}\left( {n,{\bf H}} \right)$   is  a
left linear combination of its other rows, or its $j$-th column    is  a
right linear combination of its other columns, i.e. $\exists c_{l} \in{{\bf H}}$ for $l = \overline {1,k}$ such that
${\rm {\bf a}}_{i.} = c_{1} {\rm {\bf
a}}_{i_{1} .} + \ldots + c_{k}  {\rm {\bf a}}_{i_{k} .}$ or ${\rm {\bf
a}}_{.i} = {\rm {\bf a}}_{.i_{1}}   c_{1} + \ldots + {\rm {\bf
a}}_{.i_{k}} c_{k} $  for
 $\{i,i_{l}\}\subset I_{n}$, then
 $\det{\rm {\bf A}} = 0.
$
\end{cor}
From Theorems \ref{theorem:row_combin}, \ref{theorem:colum_combin}
and basic properties of the row-column determinants for
arbitrary matrices, we can obtain the following theorems  as well.
\begin{thm} \label{theorem:row_combin_add}If the $i$-th row of a Hermitian matrix ${\rm {\bf A}}\in
{\rm M}\left( {n,{\bf H}} \right)$   is added a left
linear combination of its other rows, then \begin{gather*}
   {\rm{rdet}}_{i} \,{\rm {\bf A}}_{i \cdot } \left( {{\rm {\bf
a}}_{i.} + c_{1} \cdot {\rm {\bf a}}_{i_{1} .} + \ldots + c_{k}
\cdot {\rm {\bf a}}_{i_{k} .}} \right) =\\
  ={\rm{cdet}} _{i}
\,{\rm {\bf A}}_{i \cdot} \left( {{\rm {\bf a}}_{i.} + c_{1} \cdot
{\rm {\bf a}}_{i_{1} .} + \ldots + c_{k} \cdot {\rm {\bf
a}}_{i_{k} .}}  \right) = \det {\rm {\bf A}},
\end{gather*}
 where $ c_{l} \in {{\bf H}}$ for all $l = \overline {1,k}$ and
$\{i,i_{l}\}\subset I_{n}$.
\end{thm}
\begin{thm}\label{theorem:col_combin_add} If the $j$-th column of a Hermitian matrix ${\rm {\bf A}}\in
{\rm M}\left( {n,{\bf H}} \right)$   is added a right
linear combination of its other columns, then
\begin{gather*}
   {\rm{cdet}} _{j}\, {\rm {\bf A}}_{.j} \left( {{\rm {\bf a}}_{.j} +
{\rm {\bf a}}_{.j_{1}}
 c_{1} + \ldots + {\rm {\bf a}}_{.j_{k}}   c_{k}}
\right) = \\
  ={\rm{rdet}} _{j} \,{\rm {\bf A}}_{.j} \left( {{\rm {\bf
a}}_{.j} + {\rm {\bf a}}_{.j_{1}}   c_{1} + \ldots + {\rm {\bf
a}}_{.j_{k}}  c_{k}}  \right) = \det {\rm {\bf A}},
\end{gather*}
where $c_{l} \in{{\bf H}}$ for all $l = \overline
{1,k}$ and $\{j,j_{l}\}\subset J_{n}$.
\end{thm}

\section{ Determinantal representations the inverse of a Hermitian matrix}
\subsection{The inverse of a Hermitian matrix}
\begin{thm} \label{thm:inver_her} If  ${\rm {\bf A}}\in
{\rm M}\left( {n,{\bf H}} \right)$ is Hermitian and $\det{\bf A}\neq 0$, then there exist an unique right
inverse  matrix $(R{\rm {\bf A}})^{ - 1}$ and an unique left
inverse matrix $(L{\rm {\bf A}})^{ - 1}$ of  ${\rm {\bf A}}$, where $\left( {R{\rm {\bf A}}} \right)^{
- 1} = \left( {L{\rm {\bf A}}} \right)^{ - 1} = :{\rm {\bf A}}^{ -
1}$, and they have the following  determinantal representations, respectively,
\begin{equation}
\label{eq:repr_righ_inv_her}
   \left( {R{\rm {\bf A}}} \right)^{ - 1} = {\frac{{1}}{{\det
{\rm {\bf A}}}}}
\begin{pmatrix}
  R_{11} & R_{21} & \cdots & R_{n1}\\
  R_{12} & R_{22} & \cdots & R_{n2}\\
  \cdots & \cdots & \cdots& \cdots\\
  R_{1n} & R_{2n} & \cdots & R_{nn}
\end{pmatrix},
\end{equation}
\begin{equation}
\label{eq:repr_left_inv_her}
 \left( {L{\rm {\bf A}}} \right)^{ - 1} = {\frac{{1}}{{\det {\rm
{\bf A}}}}}
\begin{pmatrix}
  L_{11} & L_{21} & \cdots & L_{n1} \\
  L_{12} & L_{22} & \cdots & L_{n2} \\
  \cdots & \cdots & \cdots & \cdots \\
  L_{1n} & L_{2n} & \cdots & L_{nn}
\end{pmatrix},
\end{equation}
 where $R_{ij}$ and  $L_{ij}$ can be obtained by (\ref{eq:Rij}) and (\ref{eq:Lij}), respectively,  for all $ i,j =
\overline {1,n}$.
\end{thm}
\emph{Proof.} Let ${\rm {\bf B}} = {\rm {\bf A}} \cdot \left( {R{\rm
{\bf A}}} \right)^{ - 1}$.  We obtain the entries of ${\rm {\bf
B}}$ by multiplying matrices. For all $i = \overline {1,n}$, we have
\[
   b_{i{\kern 1pt} i} = \left( {\det {\rm {\bf A}}} \right)^{ -
1}{\sum\limits_{j = 1}^{n} {a_{i{\kern 1pt} j} \cdot R_{i{\kern
1pt} j}} } = \left( {\det {\rm {\bf A}}} \right)^{ -
1}{\rm{rdet}}_{i}\, {\rm {\bf A}} = {\frac{{\det {\rm {\bf
A}}}}{{\det {\rm {\bf A}}}}} = 1,\] and for all ${i \ne j}$
\[
b_{i{\kern 1pt} j} = \left( {\det {\rm {\bf A}}} \right)^{ -
1}{\sum\limits_{s = 1}^{n} {a_{i{\kern 1pt} s} \cdot R_{j{\kern
1pt} s}} } = \left( {\det {\rm {\bf A}}} \right)^{ - 1}{\rm{rdet}}
_{j} {\rm {\bf A}}_{j{\kern 1pt}.} \left( {{\rm {\bf a}}_{i{\kern
1pt}.}} \right).
\]
If $i \ne j$, then by Theorem \ref{repl_row} ${\rm{rdet}} _{j}
{\rm {\bf A}}_{j{\kern 1pt}.} \left( {{\rm {\bf a}}_{i{\kern
1pt}.}} \right)=0$. Consequently $b_{i{\kern 1pt} j} = 0$. Thus
${\rm {\bf B}} = {\rm {\bf I}}$ and $\left( {R{\rm {\bf A}}}
\right)^{ - 1}$ is the right inverse of the Hermitian matrix ${\rm
{\bf A}}$.

Suppose ${\rm {\bf D}} = \left( {L{\rm {\bf A}}} \right)^{ -
1}{\rm {\bf A}}$. We again  get the entries of ${\rm {\bf D}}$ by
multiplying matrices. For all $i = \overline {1,n}$,
\[
 d_{i{\kern 1pt} i} = \left( {\det {\rm {\bf A}}} \right)^{ -
1}{\sum\limits_{i = 1}^{n} {L_{i{\kern 1pt} j} \cdot a_{i{\kern
1pt} j}} } = \left( {\det {\rm {\bf A}}} \right)^{ - 1}{\rm{cdet}}
_{j} {\rm {\bf A}} = {\frac{{\det {\rm {\bf A}}}}{{\det {\rm {\bf
A}}}}} = 1,\] and for all $i \ne j$,
\[
  d_{i{\kern 1pt} j} = \left( {\det {\rm {\bf A}}} \right)^{ -
1}{\sum\limits_{s = 1}^{n} {L_{s{\kern 1pt} {\kern 1pt} i} \cdot
a_{s{\kern 1pt} j}} }  = \left( {\det {\rm {\bf A}}} \right)^{ -
1}{\rm{cdet}} _{i} {\rm {\bf A}}_{.i} \left( {{\rm {\bf a}}_{.j}}
\right).
\]
If $i \ne j$, then by Theorem \ref{repl_col} ${\rm{cdet}} _{i}
{\rm {\bf A}}_{.i} \left( {{\rm {\bf a}}_{.j}} \right)=0$.
Therefore $d_{i{\kern 1pt} j} = 0$ for all $i \ne j$. Thus ${\rm
{\bf D}} = {\rm {\bf I}}$ and $\left( {L{\rm {\bf A}}} \right)^{ -
1}$ is the left inverse
 of the Hermitian matrix ${\rm {\bf A}}$.

The equality $\left( {R{\rm {\bf A}}} \right)^{ - 1} = \left(
{L{\rm {\bf A}}} \right)^{ - 1}$ because of the uniqueness of inverses over associative rings.
$\Box$

Moreover, the following  criterion of invertibility of a Hermitian matrix
can be obtained.
\begin{thm}\label{theorem:inver_equiv} If ${\rm {\bf A}}\in {\rm M}\left( {n,{\rm
{\bf{H}}}} \right)$ is Hermitian, then
 the following propositions are equivalent.
\begin{itemize}
  \item[i)] ${\rm {\bf A}}$ is invertibility, i.e. ${\rm {\bf A}}\in
GL\left( {n,{\bf{H}}} \right);$
  \item[ii)]  rows of ${\rm {\bf A}}$ are left-linearly independent;
  \item[iii)]  columns of ${\rm {\bf A}}$ are right-linearly
independent;
  \item[iiii)] $\det {\rm {\bf A}}\ne 0$.
\end{itemize}
\end{thm}
\emph{Proof.} $i) \Rightarrow ii)$ Consider a right system of  linear equations ${\rm {\bf A}} \cdot {\rm {\bf x}} = {\rm
{\bf y}}$. The fact ${\rm {\bf A}}$ is invertible means that  the linear transformation ${\rm {\bf A}}: {\bf x}\rightarrow {\bf y}$ is a  bijection. Suppose that  rows of ${\rm {\bf A}}$ are left-linearly dependent. It means that  $\exists i\in I_{n}$ and $\exists c_{l} \in{{\bf H}}$ for $l = \overline {1,k}$ such that
${\rm {\bf a}}_{i.} = c_{1} {\rm {\bf
a}}_{i_{1} .} + \ldots + c_{k}  {\rm {\bf a}}_{i_{k} .}$. Then, by elementary row operations  the $i$-th  row reduce to zero, and we lose the bijectivity of the linear transformation ${\rm {\bf A}}$. It follows that ${\rm {\bf A}}$ is non-invertible. Hence, the supposition is false, and  rows of ${\rm {\bf A}}$ are left-linearly independent.

The equivalence $ii) \Rightarrow iii)$ can be proved similarly by considering a left system of  linear equations ${\rm {\bf x}}{\rm {\bf A}} = {\rm
{\bf y}}$.

The equivalences $ii) \Rightarrow iiii)$ and $iii) \Rightarrow iiii)$ follow from Corollary \ref{cor:colum_combin}.

Finally, the equivalence $iiii) \Rightarrow i)$ is given by Theorem \ref{thm:inver_her}.
$\Box$
\begin{rem} By Theorems \ref{theorem:inver_equiv},
\ref{theorem:row_combin_add} and \ref{theorem:col_combin_add}, the determinant of a coquaternionic Hermitian matrix satisfy Axioms 1,3 of a noncommutative determinant.
\end{rem}

\subsection{ Cramer's rule for  systems of linear coquaternionic equations  in Hermitian case}
\begin{thm} \label{theorem:right_system}Let
\begin{equation}
\label{eq:right_syst} {\rm {\bf A}} \cdot {\rm {\bf x}} = {\rm
{\bf y}}
\end{equation}
be a right system of  linear equations  with a matrix of
coefficients ${\rm {\bf A}}\in {\rm M}(n,{{\bf H}})$, a
column of constants ${\rm {\bf y}} = \left( {y_{1} ,\ldots ,y_{n}
} \right)^{T}\in {\rm {\bf H}}^{n\times 1}$, and a
column of unknowns ${\rm {\bf x}} = \left( {x_{1} ,\ldots ,x_{n}}
\right)^{T}$. If ${\rm {\bf A}}$ is Hermitian and ${\det}{ \rm{\bf A}} \ne 0$, then the
solution of (\ref{eq:right_syst}) is given by
components,
\begin{equation}
\label{eq:cr_right_syst_her}
x_{j} = {\frac{{{\rm{cdet}} _{j} {\rm {\bf A}}_{.j} \left( {{\rm
{\bf y}}} \right)}}{{\det {\rm {\bf A}}}}}, \quad  { j = \overline
{1,n}}.
\end{equation}
\end{thm}
\emph{Proof.} Since ${\det}{ \rm{\bf A}} \ne 0$, then, by Theorem \ref{thm:inver_her}, there exists  the unique inverse
matrix ${\rm {\bf A}}^{ - 1}$. From this the existence and
uniqueness of solutions of (\ref{eq:right_syst}) follows
immediately.

By considering  ${\rm {\bf A}}^{ - 1}$ as the left inverse, the solution of
(\ref{eq:right_syst}), ${\rm {\bf x}} = {\rm {\bf A}}^{ - 1} \cdot {\rm {\bf y}}$, can be represented  by components as follows,
\[
x_{j} = \left( {{\det}{ \rm{\bf A}}} \right)^{ -
1}{\sum\limits_{i = 1}^{n} {L_{ij} \cdot y_{i}} },\, \, j =
\overline {1,n},
\]
where $L_{ij} $ is the left $ij$-th cofactor of  $  {\bf A} $. From
here (\ref{eq:cr_right_syst_her})
 follows immediately.$\Box$

\begin{thm} Let
\begin{equation}\label{eq:left_syst}
 {\rm {\bf x}} \cdot {\rm {\bf A}} = {\rm {\bf y}}
\end{equation}
be a left  system of  linear equations
 with a
matrix of coefficients ${\rm {\bf A}}\in {\rm M}(n,{{\bf H}})$, a row of constants ${\rm {\bf y}} = \left(
{y_{1} ,\ldots ,y_{n} } \right)\in {\rm {\bf H}}^{1\times n}$, and a row of unknowns ${\rm {\bf x}}
= \left( {x_{1} ,\ldots ,x_{n}} \right)$. If ${\rm {\bf A}}$ is Hermitian and ${\det}{ \rm{\bf A}} \ne 0$, then the
solution of (\ref{eq:left_syst}) is given by
components,
\[
x_{i}^{} = {\frac{{{\rm{rdet}}_{i} {\rm {\bf A}}_{i.} \left( {{\rm
{\bf y}}} \right)}}{{\det {\rm {\bf A}}}}}, \quad  { i = \overline
{1,n}}.
\]
\end{thm}
\emph{Proof.} The proof is similar to the proof of Theorem \ref{theorem:right_system} by using (\ref{eq:repr_righ_inv_her}) for determinantal representation of ${\bf A}^{-1}$.
$\Box$
\begin{ex}
Let  consider a right system of  linear equations
\begin{equation}\label{ex:Ax}
{\bf A}{\bf x}={\bf b}
\end{equation}
with the matrix   ${\bf A}$ from (\ref{eq:A}) and  ${\bf b}=(i\,\,j\,\,k)^{T}$.
Since ${\rm {\bf A}}$ is Hermitian and ${\det}{ \rm{\bf A}}=4$, we can find the  solution of (\ref{ex:Ax}) by Cramer's rule (\ref{eq:cr_right_syst_her}).
\begin{multline*}
{x}_{1}= \frac{1}{{\det}{ \rm{\bf A}}}\,{\rm{cdet}} _{1}\begin{pmatrix}
  i & 1-k& 1-j \\
  j& 0  & 1+j\\
  k& 1-j  & 0
\end{pmatrix}=\frac{-3-i+3j+k}{4},\\
{x}_{2}= \frac{1}{{\det}{ \rm{\bf A}}}\,{\rm{cdet}} _{2}\begin{pmatrix}
  0 & i& 1-j \\
  1+k& j  & 1+j\\
  1+j& k  & 0
\end{pmatrix}=\frac{1+3i+j-k}{4},\\
{x}_{3}= \frac{1}{{\det}{ \rm{\bf A}}}\,{\rm{cdet}} _{3}\begin{pmatrix}
  0 & 1-k& i \\
  1+k& 0  & j\\
  1+j& 1-j  & k
\end{pmatrix}=\frac{2j+2k}{4}.
\end{multline*}
Now, we shall find the inverse ${\bf A}^{-1}$ of ${\bf A}$ by (\ref{eq:repr_left_inv_her}).
\begin{multline*}
L_{11}={\rm{cdet}} _{1}\begin{pmatrix}
  0 & 1+j \\
  1-j& 0
\end{pmatrix}= 0,  L_{12}=-{\rm{cdet}} _{1}\begin{pmatrix}
  1+k & 1+j \\
  1+j& 0
\end{pmatrix}= 2+2j, \\
 L_{13}=-{\rm{cdet}} _{1}\begin{pmatrix}
  1+j & 1-j \\
  1+k& 0
\end{pmatrix}= 1+i-j+k, \\  L_{21}=-{\rm{cdet}} _{1}\begin{pmatrix}
  1-k & 1-j \\
  1-j& 0
\end{pmatrix}= 2-2j,
L_{22}={\rm{cdet}} _{1}\begin{pmatrix}
  0 & 1-j \\
  1+j& 0
\end{pmatrix}= 0, \\ L_{23}=-{\rm{cdet}} _{2}\begin{pmatrix}
  0 & 1-k \\
  1+j&  1-j
\end{pmatrix}= 1+i+j-k, \\
 L_{31}=-{\rm{cdet}} _{2}\begin{pmatrix}
  0 & 1+j \\
  1-k& 1-j
\end{pmatrix}= 1-i+j-k, \\  L_{32}=-{\rm{cdet}} _{2}\begin{pmatrix}
  0 & 1-j \\
  1+k& 1+j
\end{pmatrix}= 1-i-j+k, \\
 L_{33}={\rm{cdet}} _{2}\begin{pmatrix}
  0 & 1-k \\
  1+k& 0
\end{pmatrix}= 0.
 \end{multline*}
Therefore,
\[{\bf A}^{ - 1} = {\frac{{1}}{4}}
\begin{pmatrix}
  0 & 2-2j & 1-i+j-k\\
  2+2j & 0 &  1-i-j+k\\
   1+i-j+k & 1+i+j-k & 0
\end{pmatrix}.\]
Finally, we see that by  the matrix method  the identical result is obtained,
\begin{multline*}\begin{pmatrix}
 x_{1}\\
  x_{2}\\
   x_{3}
\end{pmatrix}={\frac{{1}}{4}}
\begin{pmatrix}
  0 & 2-2j & 1-i+j-k\\
  2+2j & 0 &  1-i-j+k\\
   1+i-j+k & 1+i+j-k & 0
\end{pmatrix}\begin{pmatrix}
 i\\
  j\\
   k
\end{pmatrix}=\\{\frac{{1}}{4}}\begin{pmatrix}
 -3-i+3j+k\\
  1+3i+j-k\\
   2j+2k
\end{pmatrix}.\end{multline*}
\end{ex}

\section{Cramer's rules for some coqaternionic matrix equations}
\begin{thm} Suppose
\begin{equation}\label{equat:AXB}
{\rm {\bf A}}{\rm {\bf X}}{\rm {\bf B}} = {\rm {\bf C}}
\end{equation}
\noindent is a two-sided matrix equation, where $ {\rm
{\bf A}} \in{\bf{H}}^{m\times m}$, ${\rm {\bf B}}\in{\bf{H}}^{n\times n}$, ${\rm {\bf C}}\in{\bf{H}}^{m\times n}$ are given, ${\rm {\bf X}} \in{\bf{H}}^{m\times n}$ is unknown, and ${\bf A}$, ${\bf B}$ are Hermitian. If ${\det} {\rm {\bf A}} \ne 0$
and ${\det} {\rm {\bf B}} \ne 0$ , then the unique solution of (\ref{equat:AXB}) can be represented as follows,
\begin{equation}
\label{eq:rdet_AXB} x_{ij} = {\frac{{{\rm rdet} _{j} {\rm {\bf
B}}_{j.\,} \left( {{\rm {\bf c}}_{i\,.}^{{\rm
{\bf A}}}} \right)}}{{{\det}  {\rm {\bf A}}\cdot {\det}
{\rm {\bf B}}}}},
\end{equation}
\noindent or
\begin{equation}
\label{eq:cdet_AXB} x_{ij} = {\frac{{{\rm cdet} _{i} {\rm {\bf
A}}_{.\,i\,} \left( {{\rm {\bf c}}_{.j}^{{\rm
{\bf B}}}} \right)}}{{{\det} {\rm {\bf A}}\cdot {\det}
{\rm {\bf B}}}}},
\end{equation}
\noindent where ${\rm {\bf c}}_{i\,.}^{{\rm {\bf A}}} : = \left(
{{\rm cdet} _{i}  {\bf A}_{.\,i}
\left( {{{\rm {\bf c}}}_{.1}} \right),\ldots ,{\rm cdet}
_{i}  {\rm {\bf A}}_{.\,i} \left(
{{{\rm {\bf c}}}_{.n}} \right)} \right)\in{\bf{H}}^{n\times 1}$ is the row-vector
and  ${\rm {\bf c}}_{.j}^{{\rm {\bf B}}} : = \left( {{\rm rdet}
_{j} {\rm {\bf B}}_{j.} \left( {{{\rm
{\bf c}}}_{1.}} \right),\ldots ,{\rm rdet} _{j} {\rm {\bf
B}}_{j.} \left( {{{\rm {\bf
c}}}_{m.}^{}} \right)} \right)^{T}\in{\bf{H}}^{1\times m}$ is the column-vector and
${{\rm {\bf c}}}_{i\,.}$, ${{\rm {\bf c}}}_{.\,j}$ are
the  i-th row  and the j-th column  of ${\rm {\bf
{C}}}$, respectively, for all $i = \overline{1,m} $, $j = \overline{1,n} $.
\end{thm}
\emph{Proof.} By Theorem \ref{thm:inver_her},  ${\rm {\bf A}}$ and ${\rm {\bf B}}$ are invertible. There
exists the unique solution of
(\ref{equat:AXB}), ${\rm {\bf X}} = {\rm {\bf
A}}^{ - 1}{\rm {\bf C}}{\rm {\bf B}}^{ - 1}$. If we represent
${\rm {\bf A}}^{ - 1}$ as a left inverse by (\ref{eq:repr_left_inv_her}) and $\left( {{\rm {\bf
B}}} \right)^{ - 1}$ as a right inverse by (\ref{eq:repr_righ_inv_her}), then  we have
\begin{gather*}
{\rm {\bf X}}
 = \begin{pmatrix}
   x_{11} & x_{12} & \ldots & x_{1n} \\
   x_{21} & x_{22} & \ldots & x_{2n} \\
   \ldots & \ldots & \ldots & \ldots \\
   x_{m1} & x_{m2} & \ldots & x_{mn} \
 \end{pmatrix}
 = {\frac{{1}}{{{\det} {\rm {\bf A}}}}}\begin{pmatrix}
  {L} _{11}^{{\rm {\bf A}}} & {L} _{21}^{{\rm {\bf A}}}&
   \ldots & {L} _{m1}^{{\rm {\bf A}}} \\
  {L} _{12}^{{\rm {\bf A}}} & {L} _{22}^{{\rm {\bf A}}} &
  \ldots & {L} _{m2}^{{\rm {\bf A}}} \\
  \ldots & \ldots & \ldots & \ldots \\
 {L} _{1m}^{{\rm {\bf A}}} & {L} _{2m}^{{\rm {\bf A}}}
  & \ldots & {L} _{mm}^{{\rm {\bf A}}}
\end{pmatrix}\times\\
  \times\begin{pmatrix}
    {c}_{11} & {c}_{12} & \ldots & {c}_{1n} \\
    {c}_{21} & {c}_{22} & \ldots & {c}_{2n} \\
    \ldots & \ldots & \ldots & \ldots \\
    {c}_{m1} & {c}_{m2} & \ldots & {c}_{mn} \
  \end{pmatrix}
{\frac{{1}}{{{\det} {\rm {\bf B}}}}}
\begin{pmatrix}
 {R} _{\, 11}^{{\rm {\bf B}}} & {R} _{\, 21}^{{\rm {\bf B}}}
 &\ldots & {R} _{\, n1}^{{\rm {\bf B}}} \\
 {R} _{\, 12}^{{\rm {\bf B}}} & {R} _{\, 22}^{{\rm {\bf B}}} &\ldots &
 {R} _{\, n2}^{{\rm {\bf B}}} \\
 \ldots  & \ldots & \ldots & \ldots \\
 {R} _{\, 1n}^{{\rm {\bf B}}} & {R} _{\, 2n}^{{\rm {\bf B}}} &
 \ldots & {R} _{\, nn}^{{\rm {\bf B}}}
\end{pmatrix},
\end{gather*}
 where ${L}_{ij}^{{\rm {\bf A}}} $ is a left $ij$-th
 cofactor of ${\rm {\bf A}}$ for all $i,j = \overline{1,m} $, and ${R}_{i{\kern 1pt} j}^{{\rm {\bf
B}}} $ is a right  $ij$-th cofactor of ${\rm {\bf B}}$ for all $i,j = \overline{1,n} $. It implies
\begin{equation}
\label{eq:sum_AXB} x_{ij} = {\frac{{{\sum\limits_{l = 1}^{n}
{\left( {{\sum\limits_{k = 1}^{m} {{L}_{ki}^{{\rm {\bf A}}}
{c}_{kl}} } } \right)}} {R}_{jl}^{{\rm {\bf B}}}} }{{{\det} {\rm {\bf A}}\cdot{\det} {\rm {\bf B}}}}},
\end{equation}
 for all $i = \overline{1,m} $, $j = \overline{1,n} $. From this by Definition
\ref{def_cdet}, we obtain
\[{\sum\limits_{k = 1}^{m}
{{L}_{ki}^{{\rm {\bf A}}} {c}_{k\,l}} }  = {\rm cdet} _{i}
{\rm {\bf A}}_{.\,i} \left( {{{\rm
{\bf c}}}_{.l}} \right),\]
\noindent where ${\rm {\bf
c}}_{.\,l} $ is the $l$-th column of ${{\rm {\bf
C}}}$ for all $l = \overline{1,n} $. Consider the row-vector
 \[{\rm {\bf
c}}_{i\,.}^{{\rm {\bf A}}} : = \left( { {\rm cdet} _{i} {\rm {\bf
A}}_{.\,i\,} \left( {{{\rm {\bf
c}}}_{.1}} \right),\ldots , {\rm cdet} _{i} {\rm {\bf A}}_{.\,i} \left( {{{\rm {\bf c}}}_{.n}} \right)}
\right)\]
  for all $i = \overline{1,m} $. By  Definition
\ref{def_rdet},
${\sum\limits_{l = 1}^{n} c_{il}^{{\bf A}}
{R}_{jl}^{{\rm {\bf B}}}}={{\rm rdet} _{j} {\rm {\bf
B}}_{j.} \left( {{\rm {\bf c}}_{i\,.}^{{\rm
{\bf A}}}} \right)} $, then we get
(\ref{eq:rdet_AXB}).
Having changed the order of summation in (\ref{eq:sum_AXB}), we
have
\[
x_{ij} = {\frac{{{\sum\limits_{k = 1}^{m} {{ L}_{ki}^{{\rm {\bf
A}}}} }\left( {{\sum\limits_{l = 1}^{n} {{c}_{\,kl}} }
{R}_{jl}^{{\rm {\bf B}}}} \right)}}{{{\det} {\rm {\bf
A}}\cdot{\det}{\rm {\bf B}}}}}.\]
By   Definition
\ref{def_rdet}, we have ${\sum\limits_{l = 1}^{n}
{c_{k\,l}} } {R}_{j\,l}^{{\rm {\bf B}}} = {\rm rdet} _{j} {\rm
{\bf B}}_{j.} \left( {{{\rm {\bf
c}}}_{k.}} \right)$, where ${{\rm {\bf c}}}_{k.} $ is the
$k$-th row-vector of ${{\rm {\bf C}}}$ for all $k = \overline{1,n}
$. Denote  the following column-vector by
\[{\rm {\bf c}}_{.j}^{{\rm {\bf B}}} : = \left(
{{\rm rdet} _{j} {\rm {\bf B}}_{j.} \left(
{{{\rm {\bf c}}}_{1.}} \right),\ldots ,{\rm rdet} _{j}
{\rm {\bf B}}_{j.\,} \left( {{{\rm {\bf
c}}}_{m.}} \right)} \right)^{T}\]
  for all
$j = \overline{1,n} $. By  Definition \ref{def_cdet}, ${\sum\limits_{k = 1}^{n} {{
L}_{ki}^{{\rm {\bf A}}}} } c_{kj}^{\bf B}={{{\rm cdet} _{i} {\rm {\bf
A}}_{.\,i\,} \left( {{\rm {\bf c}}_{.j}^{{\rm
{\bf B}}}} \right)}}$, then we finally have
(\ref{eq:cdet_AXB}).
$\Box$

If,  in (\ref{equat:AXB}), we put $ {\bf
A}={\bf I}_m$ or  $ {\bf B}={\bf I}_n$, then, respectively,  we evidently get the following  corollaries.
\begin{cor} Suppose
\begin{equation}\label{eq:AX}
 {\rm {\bf A}}{\rm {\bf X}} = {\rm {\bf C}}
\end{equation}
\noindent is a right matrix equation, where $ {\bf
A} \in {\bf{H}}^{m\times m}$, ${\bf C}\in {\bf{H}}^{m\times n}$
are given, ${\rm {\bf X}} \in {\bf{H}}^{m\times n}$ is
unknown, and ${\bf A}$ is Hermitian. If ${\det} {\rm {\bf A}} \ne 0$, then the unique solution of (\ref{eq:AX}) can be represented as follows,
\begin{equation*}
 x_{i\,j} = {\frac{{{\rm cdet} _{i} {\rm {\bf
A}}_{.i} \left( {{{\rm {\bf c}}}_{.j}}
\right)}}{{ {\det} {\rm {\bf A}}}}}
\end{equation*}
\noindent where ${{\rm {\bf c}}}_{.j}$ is the $j$-th column of
${{\rm {\bf C}}}$, for all $i=\overline{ 1,m}$, $j =\overline{1,n}$.
\end{cor}
\begin{cor} Suppose
\begin{equation}\label{eq:XB}
 {\rm {\bf X}}{\rm {\bf B}} = {\rm {\bf C}}
\end{equation}
 is a left matrix equation, where ${\rm {\bf B}}\in {\bf{H}}^{n\times n}$, ${\rm {\bf C}}\in {\bf{H}}^{m\times n}$
are given, ${\rm {\bf X}} \in {\bf{H}}^{m\times n}$ is
unknown, and  ${\bf B}$ is Hermitian. If ${\det} {\rm {\bf B}} \ne 0$, then the unique solution of (\ref{eq:XB}) can be represented as follows,
\begin{equation*}
 x_{ij} = {\frac{{{\rm rdet} _{j} {\rm {\bf
B}}_{j.\,} \left( {{{\rm {\bf
c}}}_{i\,.}} \right)}}{{{\det} {\rm {\bf B}}}}}
\end{equation*}
\noindent where ${{\rm {\bf c}}}_{i.}$ is the $i$-th row
of ${{\rm {\bf C}}}$, for all $i=\overline{ 1,m}$, $j =\overline{1,n}$.
\end{cor}
\begin{ex}
Let  consider the matrix equations
\begin{equation}\label{ex:AXB}
{\rm {\bf A}}{\rm {\bf X}}{\rm {\bf B}} = {\rm {\bf C}}
\end{equation}
with the matrix   ${\bf A}$ from (\ref{eq:A}),   ${\bf B}=\begin{pmatrix}
1& k\\
-k & 1\end{pmatrix}$, and ${\bf C}=\begin{pmatrix}
i& 1\\
0 & j\\
k& -i\end{pmatrix}.$ Since ${\rm {\bf A}}$, ${\rm {\bf B}}$ are Hermitian and ${\det}{ \rm{\bf A}}=4$, and ${\det}{ \rm{\bf B}}=2$, we can find the  solution of (\ref{ex:AXB}) by Cramer's rule (\ref{eq:cdet_AXB}). Firstly, we obtain the column-vectors ${\rm {\bf c}}_{.j}^{{\rm {\bf B}}}$ for $j=\overline{1,2}$. Since
\begin{multline*}
c_{11}^{\bf B}={\rm{rdet}} _{1}{\bf B}_{.1}({\bf c}_{.1})={\rm{rdet}} _{1}\begin{pmatrix}
  i & 1 \\
  -k& 1
\end{pmatrix}= i+k,\\  c_{21}^{\bf B}={\rm{rdet}} _{1}{\bf B}_{.1}({\bf c}_{.2})={\rm{rdet}} _{1}\begin{pmatrix}
  0 & j \\
  -k& 1
\end{pmatrix}= -i, \\
c_{31}^{\bf B}={\rm{rdet}} _{1}{\bf B}_{.1}({\bf c}_{.3})={\rm{rdet}} _{1}\begin{pmatrix}
  k & -i \\
  -k& 1
\end{pmatrix}= j+k,
 \end{multline*}
 then ${\rm {\bf c}}_{.1}^{{\rm {\bf B}}}=\begin{pmatrix}
  i+k \\
  -i\\
  j+k
\end{pmatrix}$. Similarly, we get
 ${\rm {\bf c}}_{.2}^{{\rm {\bf B}}}=\begin{pmatrix}
  1-j \\
  j\\
  -1-i
\end{pmatrix}$. Then by (\ref{eq:cdet_AXB}), we have
\begin{equation*}
x_{11}= {\frac{{{\rm cdet} _{1} {\rm {\bf
A}}_{.\,1} \left( {{\rm {\bf c}}_{.1}^{{\rm
{\bf B}}}} \right)}}{{{\det} {\rm {\bf A}}\cdot {\det}
{\rm {\bf B}}}}}=\frac{1}{8}\,{\rm{cdet}} _{1}\begin{pmatrix}
  i+k & 1-k & 1-j\\
  -i& 0 & 1+j\\
  j+k &1-j  & 0
  \end{pmatrix}= \frac{-2i+j-k}{4}.\end{equation*}
 Similarly, we obtain
 \begin{multline*}
 x_{12}= \frac{-2+j+k}{4}, x_{21}= \frac{i+j}{4}, x_{22}= \frac{-1-k}{4},
 x_{31}= \frac{1+i+j+3k}{8},\\ x_{32}= \frac{3-2i-j+2k}{8}.
 \end{multline*}
 Finally, \[{\bf X}=\frac{1}{8}\begin{pmatrix}
  -4i+2j-2k &-4+2j+2k \\
  2i+2j & -2-2k\\
  1+i+j+3k & 3-2i-j+2k
\end{pmatrix}.\]
\end{ex}


\begin{thebibliography}{9}

\bibitem{le}   D.W. Lewis,  \textit{Quaternion algebras and the algebraic legacy
of Hamilton's quaternions}. Irish Math. Soc. Bulletin \textbf{ 57} (2006),
 41-64.
\bibitem{coc}J. Cockle,  \textit{On systems of algebra involving more than one imaginary}. Phil. Mag. \textbf{35} (1849), 434-435.
\bibitem{ku}
L. Kula, Y. Yayli, \textit{Split quaternions and rotations in semi
euclidean space $E^{4}_{2}$}. J. Korean Math. Soc. \textbf{ 44} (2007),  1313-1327.

\bibitem{oz1}M.\"{O}zdemir, A.A. Ergin, \textit{Rotations with unit timelike quaternions in Minkowski 3-space.} J. Geom.Phys. \textbf{56} (2006), 322–336.
\bibitem{br}
P. Bracken, \textit{Split-quaternionic representation of the moving frame for
timelike surfaces in 3-dimensional Minkowski spacetime}. Journal of Mathematics and Statistics \textbf{6}(1) (2010), 56-59.
\bibitem{sh} V. Shpakivsky, \textit{Linear quaternionic equations and their systems}. Adv. Appl. Clifford Algebras \textbf{21} (2011) 637-645.
\bibitem{er}M. Erdo\u{g}du, M. \"{O}zdemir, \textit{Two-sided linear split quaternionic
equations with unknowns}. Linear and Multilinear Algebra \textbf{63} (2015), 97-106.


 \bibitem{al}
Y. Alag\"{o}z, K. H. Oral, S.Y\"{u}ce, \textit{Split quaternion matrices}. Miskolc Mathematical Notes \textbf{  13}  (2012), 223-232.

 \bibitem{ky}I. Kyrchei, \textit{The column and row immanants of matrices over a split quaternion algebra}. Adv. Appl. Clifford Algebras \textbf{25} (2015), 611-619.



\bibitem{ky1} I. Kyrchei, \textit{Cramer's rule for quaternionic systems of linear equations}. Fundamentalnaya i Prikladnaya Matematika \textbf{13}4 (2007), 67-94.
 \bibitem{ky2}     I. Kyrchei, \textit{The theory of the column and row determinants in a quaternion linear algebra}.  In: Albert R. Baswell (Eds.), Advances in Mathematics Research 15,  pp. 301-359,  Nova Sci. Publ., New York, 2012.
\bibitem{ky3} I. Kyrchei, \textit{Determinantal representations of the Moore-Penrose inverse over the quaternion skew field}. Journal of Mathematical Sciences \textbf{180} (012),  23-33.

\bibitem{ky4}    I. Kyrchei, \textit{Explicit representation formulas for the minimum norm least squares solutions of some quaternion matrix equations}.  Linear Algebra Appl. \textbf{438} (2013),  136-152.
   \bibitem{ky5} I. Kyrchei,  \textit{Determinantal representations of the Drazin inverse over the quaternion skew field with applications to some matrix equations}. Appl. Math. Comp. \textbf{238} (2014), 193-207.
  \bibitem{ky6} I. Kyrchei,  \textit{Determinantal representations of the W-weighted Drazin inverse over the quaternion skew field}. Appl. Math. Comp. \textbf{264} (2015), 453-465.

\bibitem{so1} G.J. Song, Q.W. Wang, H.X. Chang, \textit{Cramer rule for the unique
solution of restricted matrix equations over the quaternion skew
field}. Comput. Math. Appl. \textbf{61}  (2011), 1576-1589.

\bibitem{so2} G.J. Song, Q.W. Wang, \textit{Condensed Cramer rule for some restricted
quaternion linear equations}. Appl. Math. Comp. \textbf{218} (2011),
 3110-3121.

\bibitem{so3} G.J. Song, \textit{Bott-Duffin inverse over the quaternion skew field with applications}.  Journal of Applied Mathematics and Computing \textbf{41} (2013),  377-392.


\bibitem{er1}
 M. Erdo\"{g}du,
M. O\"{z}demir,  \textit{On complex split quaternion matrices}.
Adv. Appl. Clifford Algebras
\textbf{23}  (2013), 625-638.
\bibitem{er2}
 M. Erdo\"{g}du,
M. O\"{z}demir,
\textit{On eigenvalues of split quaternion matrices}.
Adv. Appl. Clifford Algebras
\textbf{ 23}  (2013),  615-623.


 \bibitem{cl}  C. Flaut, V. Shpakivskyi, \textit{On complex split quaternion matrices}.
Adv. Appl. Clifford Algebras
\textbf{23} ( 2013),  657-671.

\bibitem{as} H.Aslaksen,  \textit{Quaternionic determinants}. Math.
Intellig. \textbf{18} (1996), 57-65.
 \bibitem{co} N.Cohen, S. De Leo, \textit{The quaternionic determinant}. Elec. J. Lin.
Alg. \textbf{7} (2000), 100-111.

\bibitem{ge1} I. Gelfand, V. Retakh, \textit{A
 determinants of matrices over noncommutative rings}. Funct. Anal. Appl. \textbf{25} (1991), 13-35.
\bibitem{ge2} I. Gelfand, V. Retakh, \textit{A
theory of noncommutative determinants and characteristic functions
of graphs}. Funct. Anal. Appl. \textbf{26} (1992),
1-20.
\bibitem{mo}
E.H. Moore, \textit{On the determinant of an Hermitian matrix of
quaternionic elements.} Bull. Amer. Math. Soc. \textbf{28} (1922), 161-162.
    \bibitem{dy} F. J.Dyson, \textit{Quaternion determinants.} Helvetica Phys.
Acta \textbf{45} (1972), 289-302.



\end{thebibliography}
\end{document}